\documentstyle[amscd,amssymb,verbatim,12pt]{amsart}

\newtheorem{thm}{Theorem}[section]
\newtheorem{cor}[thm]{Corollary}
\newtheorem{lem}[thm]{Lemma}
\newtheorem{prop}[thm]{Proposition}

\theoremstyle{remark}
\newtheorem{remark}[thm]{Remark}

\theoremstyle{remark}

\numberwithin{equation}{section}

  \newcommand {\C}{{\mathbb C}}
  \newcommand {\N}{{\mathbb N}}
  \newcommand {\R}{{\mathbb R}}
  \newcommand {\Z}{{\mathbb Z}}

\renewcommand {\H}{{\mathcal H}}

\renewcommand {\L}{{\mathcal L}}
  \newcommand {\E}{{\mathcal E}}
  \newcommand {\F}{{\mathcal F}}
  \newcommand {\K}{{\mathcal K}}

 \newcommand {\cQ}{{\mathcal Q}}
 \newcommand {\Spec}{\operatorname{Spec}}
 \newcommand {\X }{{\mathfrak X}}
 \newcommand {\Cf}{\operatorname{Conf}}

\renewcommand{\Re}{\operatorname{Re}}

\newcommand{\Tr}{\operatorname{Tr}}
\newcommand{\Td}{\operatorname{Td}}
\newcommand{\ch}{\operatorname{ch}}

\begin{document}
\setcounter{page}{1}

\title[Extremal metrics and analytic torsion]{Extremal K\"ahler 
Metrics and Ray-Singer Analytic Torsion}

\author{Werner M\"uller}
\address{Universit\"at Bonn\\
Mathematisches Institut\\
Beringstrasse 1\\
D -- 53115 Bonn, Germany}
\email{mueller@@math.uni-bonn.de}

\author{Katrin Wendland}
\address{Universit\"at Bonn\\
Physikalisches Institut\\
Nu{\ss}allee 12\\
D -- 53115 Bonn, Germany}
\email{wendland@@avz109.physik.uni-bonn.de}

\subjclass{Primary: 58G26; Secondary: 58E11, 53C55}
\date{January 15, 1999}

\begin{abstract}
Let $(X,[\omega])$ be a compact K\"ahler manifold with a fixed K\"ahler class
$[\omega]$. Let $\K_\omega$ be the set of all K\"ahler metrics on $X$ whose
K\"ahler class equals $[\omega]$. In this paper we investigate the critical
points of the functional $g\in\K_\omega\mapsto\cQ(g)=\parallel
v\parallel_gT_0(X,g)^{1/2}$,  where $v$ is a fixed nonzero vector of 
the determinant
line 
$\lambda(X)$ 
associated to $H^*(X)$ and $T_0(X,g)$  is the Ray-Singer analytic torsion. For
a polarized algebraic manifold $(X,L)$ we consider a twisted 
version $\cQ_L(g)$ of this functional and assume that $c_1(L)=[\omega]$. Then
the critical points of $\cQ_L$ are exactly the  metrics 
$g\in\K_\omega$ of 
constant 
scalar curvature. In particular, if $c_1(X)=0$ or if $c_1(X)<0$ and 
$\frac{1}{2\pi}[\omega]=-c_1(X)$, then $\K_\omega$ contains a unique 
K\"ahler-Einstein metric $g_{KE}$ and $\cQ_L$ attains its absolut maximum at
$g_{KE}$. 
\end{abstract}

\maketitle

\setcounter{section}{-1}
\section{Introduction}

Let $X$ be a closed oriented surface. Given a Riemannian metric $g$ on $X$, we
denote by $\det\Delta_g$ the zeta regularized determinant of the Laplacian
$\Delta_g$ associated to $g$. Let $g_0$ be a fixed metric on $X$ and consider
the conformal equivalence class $\text{Conf}(g_0)$ of $g_0$. Let
$\text{Conf}_0(g_0)$ 
be the subset of metrics $g\in\text{Conf}(g_0)$ with
$\text{Area}(X,g)=\text{Area}(X,g_0)$. In \cite{OPS1}, Osgood, Phillips and
Sarnak studied the functional
\begin{equation}\label{0.1}
h\colon g\in\text{Conf}_0(g_0)\mapsto\det\Delta_g\in\R.
\end{equation}
One of the main results of \cite{OPS1} states that $h$ has a unique maximum 
and
this maximum is attained at the metric $\tilde g\in\text{Conf}_0(g_0)$ of
constant Gauss curvature. As a byproduct, this leads to a new proof of
the Riemann uniformization theorem.

The present paper grew out of an attempt to generalize  this work
 of Osgood, Phillips
and Sarnak  to higher dimensions. There exist different possibilities for 
doing this. In \cite{BCY}, \cite{CY},
for example,  Branson, Chang and Yang studied the analogous problem on
 four-manifolds. More precisley, on a given four-manifold, the authors 
 consider 
natural, conformally covariant differential operators 
such as 
the conformal Laplacian $\Delta+(n-2)s_g/4(n-1)$, and investigate 
extremals of the zeta function 
determinants   of  such operators in a given  conformal equivalence class
of metrics.

We take a different point of view. First of all, the conformal equivalence
class 
$\text{Conf}(g_0)$  determines a unique
complex structure on $X$ so that $g_0$ is hermitian with respect to this complex
structure. Since $\dim_\C X=1$, $g_0$ is a K\"ahler metric. Let $[\omega_0]\in
H^{1,1}(X)$ be the K\"ahler class of $g_0$ and let $\K_{\omega_0}$ be the space
of 
all K\"ahler metrics on $X$ with K\"ahler class equal to $[\omega_0]$. Then 
$\text{Conf}_0(g_0)=\K_{\omega_0}$. So, we may regard (\ref{0.1}) as a
functional on 
$\K_{\omega_0}$. The variational formula is different from the Polyakov
formula, 
but, of course, it gives the same result as in \cite{OPS1}.

This is our starting point for higher dimensional generalizations. We
consider a compact K\"ahler manifold $X$ of dimension $n$ and we fix a
K\"ahler class $[\omega]\in H^{1,1}(X)$. As above, let $\K_\omega$ be the
space 
of all K\"ahler metrics on $X$ whose K\"ahler class is equal to
$[\omega]$. For a given K\"ahler metric $g$, let $T_0(X,g)$ be the
Ray-Singer analytic torsion associated to $g$ \cite{RS}. This is a certain
weighted product of the regularized determinants of the Dolbeault-Laplace
operators $\Delta_{0,q}$, $q=1,...,n$. If $\dim_\C X=1$, then
$T_0(X,g)=c(\det\Delta_g)^{1/2}$ for some constant $c\not=0$. Thus, we may regard the
functional  
$$\tau: g\in\K_\omega\mapsto T_0(X,g)\in\R$$
as a higher dimensional analogue of (\ref{0.1}). However since, in general,
the harmonic $(0,q)$-forms vary nontrivially for $q>0$, we need to modify this
functional appropriately. Let $\parallel\cdot\parallel_{Q,g}$ be the Quillen
metric on the determinant line
\begin{equation}\label{0.2}
\lambda=\bigotimes_{q=0}^n\left(\det H^q(X)\right)^{(-1)^{(q+1)}}
\end{equation}
\cite{BGS3}.
Fix $v\in\lambda$, $v\not=0$, and put $\cQ(X,g):=\parallel v\parallel_{Q,g}$.
Then we consider the functional
\begin{equation}\label{0.3}
\cQ\colon g\in\K_\omega\mapsto \cQ(X,g)\in\R.
\end{equation}
Now recall that any variation $g_u$, $u\in(-\varepsilon,\varepsilon)$,  of $g\in\K_\omega$ is of the form
$\omega_u=\omega+\partial\overline\partial\varphi_u$ for some $\varphi_u\in
C^\infty(X)$, $u\in(-\varepsilon,\varepsilon)$. 
Using the results of Bismut, Gillet, and Soul\'e \cite{BGS3}, in section 2  we
compute the variation  $\delta\cQ/\delta\varphi$ for $\varphi\in
C^\infty(X)$.

In section 3, we briefly discuss the case of a Riemann surface and the relation
with \cite{OPS1}.

In section 4 we assume that $X$ admits a metric of constant holomorphic
curvature. Such manifolds may be regarded as higher-dimensional analogues of
Riemann surfaces. Applying the variational formula of section 2, we show that
any metric $g_{KE}\in\K_\omega$ of constant holomorphic curvature is a
critical 
point of $\cQ$. The question, if in this case $g_{KE}$ is the only critical
point of $\cQ$, remains open. Also we do not know if $g_{KE}$ is an extremum
of $\cQ$.

In section 5 we consider $K3$ surfaces. These are examples which do not admit
any metric of constant holomorphic sectional curvature. But, on the other
hand, by Yau's theorem 
\cite{Y}, every K\"ahler class $[\omega]$ on a $K3$ surface contains a
unique 
K\"ahler-Einstein metric $g_{KE}$. This metric should be a natural candidate
for a critical point of $\cQ:\K_\omega\to\R$. It follows from the variational
formula of section 2 that $g_{KE}$ is a critical point of $\cQ$ if and only if
the Euler-Lagrange equation
\begin{equation}\label{ELE}
\Delta c_2(\Omega_{KE})=0
\end{equation}
holds. We do not know if there exists any $K3$ surface that admits a
K\"ahler-Einstein metric satisfying (\ref{ELE}). However, we can show that
there are K\"ahler-Einstein metrics on certain $K3$ surfaces which do not
satisfy (\ref{ELE}).

In section 6 we introduce a modification of our functional which has a simpler
variational formula. For this purpose we assume that $X$ is a complex
projective algebraic 
manifold. Then there exists a positive line bundle $L$ over $X$. We choose the
K\"ahler form $\omega$  such that
$[\omega]=c_1(L)$. Let $g\in\K_{\omega}$ be the  metric corresponding to
$\omega$. 
 Following Donaldson \cite{D} we introduce a certain
virtual bundle 
\begin{equation}\label{0.5a}
\E=\bigoplus^{4\kappa_2(n+1)}(L-L^{-1})^{\otimes n}\oplus\bigoplus^{-\kappa_1n}
(L-L^{-1})^{\otimes (n+1)},
\end{equation}
where $\kappa_1$ and $\kappa_2$ are integers which are defined by 
$\kappa_1=\int_Xc_1(\Omega)\wedge\omega^{n-1}$ and 
$\kappa_2=\int_X\omega^n$, respectively. Similarly to (\ref{0.2}) there is a 
determinant line
$\lambda(\E)$. 
 Moreover, any hermitian metric $h$ on $L$
determines a hermitian metric $h^\E$ on $\E$. Let $\tilde
g\in\K_{\omega}$. With 
respect to the metrics 
$(\tilde g,h^\E)$, we form the analytic torsion $T_0(X,\E,\tilde g,h^\E)$ of
$X$ with 
coefficients in $\E$ and the Quillen norm $\parallel\cdot\parallel_{Q,\tilde
  g,h^\E}$ 
on $\lambda(\E)$. Fix $v\in\lambda(\E)$, $v\not=0$, and set
$$\cQ_L(\tilde g,h):=\parallel v\parallel_{Q,\tilde g,h^\E}.$$
This functional of $(\tilde g,h)$ can be turned into a functional on
$\K_{\omega}$ 
as 
follows. According 
to our choice of $\omega$, there exists a hermitian metric $h$ on $L$ such
that the curvature $\Theta_{h}$ of $h$ satisfies $\Theta_{h}=-2\pi
i\omega$. Then $h$ is determined by $g$ up to multiplication by a positive scalar. Let
$$\H_{\omega}=\{\varphi\in C^\infty(X)\mid
\omega+i\partial\overline\partial\varphi>0\}$$
be the space of K\"ahler potentials.  Given $\varphi\in\H_\omega$, let
$\omega(\varphi)=\omega+i\partial\overline\partial\varphi$ and let
$g(\varphi)$ be the metric with K\"ahler form $\omega(\varphi)$. It is well
known that the map $\varphi\in\H_\omega\mapsto g(\varphi)\in\K_\omega$ is
surjective.  
Furthermore, let $h(\varphi)=e^{-2\pi\varphi}h$. Then we put
$$\cQ_L(\varphi)=\cQ_L(g(\varphi),h(\varphi))^{(-1)^n},
\quad\varphi\in\H_\omega.$$ 
Using the variational formula established in \cite{BGS3}, one can show that
 the functional $\cQ_L$ satisfies $\cQ(\varphi+c)=\cQ(\varphi)$ for
 all $c\in\R$ and hence, it can be pushed down to a  functional on
$\K_\omega$. By the same variational formula  one can compute the variation 
of $\cQ_L$. Let $\varphi_t$,
$t\in(a,b)$, be a smooth path in $\H_\omega$ and set $g_t=g(\varphi_t)$. Then 
\begin{equation}\label{0.4}
\frac{\partial\cQ_L}{\partial
  t}(g_t)=c_n\int_X\dot\varphi_t(s(\varphi_t)-s_0)\omega(\varphi_t)^n,
\end{equation}
where $s(\varphi_t)$ is the scalar curvature of the metric $g(\varphi_t)$,
$s_0$ is the normalized total scalar curvature and $c_n>0$ is a certain
constant that depends only on $n$ and $[\omega]$. 
It follows
from (\ref{0.4}) that the critical points of $\cQ_L$ are exactly the metrics
$g\in\K_\omega$ of constant scalar curvature. Any such metric is an extremal
metric 
in the sense of Calabi \cite{Ca}.

Furthermore, formula (\ref{0.4}) agrees, up to the constant $c_n$ and the
sign,  with the  
variation of
the K-energy   $\mu\colon \K_\omega\to\R$ introduced by Mabuchi
\cite{Ma1}. The  
K-energy has been studied by Bando and Mabuchi \cite{B}, \cite{BM}, 
 \cite{Ma2}, mainly in connection with Futaki's obstruction to the existence of
 K\"ahler-Einstein metrics in the case $c_1(X)>0$.  In \cite{Ma2}, Mabuchi
 defined a natural Riemannian structure on 
 $\K_\omega$. He proved that the  sectional curvature of $\K_\omega$ is
 nonpositive  
 and that Hess($\mu$) is positive semidefinite everywhere. 
  Therefore $\cQ_L$ is also a convex functional. But this is not sufficient to
determine the type of the critical point.

In section 7 we assume that $\K_\omega$ contains a K\"ahler-Einstein
metric. Then 
we can say more about the critical points of $\cQ_L$. The
main result is the following theorem.

\begin{thm}\label{th0.1} Let $(X,L)$ be a polarized projective
  algebraic manifold. Choose a K\"ahler form $\omega$ on $X$ such that
 $[\omega]=c_1(L)$. Assume
that $\K_\omega$ contains a K\"ahler-Einstein metric $g_{KE}$. If
$c_1(X)\le0$, then $\cQ_L$ has a unique maximum which is attained at 
$g_{KE}$. If $c_1(X)>0$, then $\cQ_L$ attains its absolute maximum on the
subset $\K_{KE}\subset\K_\omega$ of K\"ahler-Einstein metrics.
\end{thm}
To prove the first part of Theorem \ref{th0.1}, we use the  evolution of the
metric by the 
complex analogue of Hamilton's Ricci flow
equation 
\begin{equation}\label{0.5}
\frac{\partial \tilde g_{i\overline j}}{\partial t}=-\tilde r_{i\overline
  j}+\frac{s_0}{2n}\tilde 
g_{i\overline j},\quad \tilde g_{i\overline j}(0)=g_{i\overline j},
\end{equation}
where $\tilde r
_{ij}$ is the Ricci tensor of $\tilde g$ and $s_0$ the
normalized total scalar curvature.  Cao \cite{C} proved that
(\ref{0.5}) has a unique solution $\tilde g_{i\overline j}(t)$ which exists
for all time and, furthermore,  if $c_1(X)\le0$, then $\tilde
g(t)$ converges to $g_{KE}$ as $t\to\infty$. Now the main observation is that
the K-energy decreases along the flow generated by (\ref{0.5}).

The case $c_1(X)>0$ is more complicated. First of all, there are obstructions
to the existence of K\"ahler-Einstein metrics \cite{Fu}, \cite{Ti2}. 
Moreover, the subspace  $\K_E\subset \K_\omega$ of
K\"ahler-Einstein metrics may have positive dimension. If
$\K_E\not=\emptyset$, it was proved in
\cite{BM} and \cite{B} that 
$\mu$  takes its minimum on $\K_E$. This implies the second statement of the
theorem.

Theorem \ref{th0.1} has some implication 
for the spectral 
determination of K\"ahler-Einstein metrics.
As above, let $g$ be a K\"ahler metric on $X$ such that $[\omega_g]=c_1(L)$
and pick a hermitian metric $h$ on $L$ such that $\Theta_h=-2\pi i\omega_g$. 
Recall that $h$ is uniquely determined by $g$ up to multiplication by a 
positive constant.
Let $q\in\{0,...,n\}$ and $k\in\N$. The 
Dolbeault-Laplace operator in $\Lambda^{0,q}(X,L^{\otimes k})$, associated to
$(g,h)$,  remains unchanged if we multiply $h$ by $\lambda\in\R^+$. 
Therefore, it is uniquely determined by $g$ and 
will be denoted by $\Delta^g_{0,q,k}$. Let  Spec($\Delta_{0,q,k}$) denote the
 spectrum of $\Delta^g_{0,q,k}$. By the above, the 
spaces $\H^{0,q}(X,L^{\otimes k})$ of 
$L^{\otimes k}$-valued 
harmonic $(0,q)$-forms also depend only on $g$. Although the inner product in 
$\H^{0,q}(X,L^{\otimes k})$ associated to $(g,h)$ depends on $h$, the 
induced $L^2$-norm on $\lambda(\E)$ is invariant under multiplication of $h$ 
by positive scalars and therefore, depends only on $g$. We denote it by 
$\parallel\cdot \parallel_g$. 

\begin{cor}\label{cor0.2}
Let $(X,L)$ be a polarized projective algebraic manifold and suppose that 
$c_1(X)\le0$. 
Choose a K\"ahler form $\omega$ on $X$ such that $[\omega]=c_1(L)$. Let
$g_{KE}$ be the unique K\"ahler-Einstein  metric in $\K_\omega$. Let 
$g\in\K_\omega$ and suppose that the following holds:
\begin{enumerate}
\item[1)]$\Spec(\Delta^g_{0,q,k})
=\Spec(\Delta^{g_{KE}}_{0,q,k})$ for all $q=0,...,n$ and 
$k=-(n+1),...,n+1$.
\item[2)]$\parallel\cdot\parallel_g=\parallel\cdot\parallel_{g_{KE}}$. 
\end{enumerate}
Then $g=g_{KE}$.
\end{cor} 

It seems to be likely that the corollary can be improved. One would expect that
$g_{KE}$ is already uniquely determined by the spectra of the Dolbeault-Laplace
operators.

In the final section 8 we  briefly discuss the relation to moduli spaces. 
These are slight modifications of results due to Fujiki and Schumacher 
\cite{FS}. First we consider a metrically polarized family of compact Hodge 
manifolds $(\pi:\X\to S,\tilde\omega,\L)$. Using $\L$, we introduce the 
family version of the virtual bundle (\ref{0.5a}) which we also denote by 
$\E$. Associated to $\E$  is the determinant line bundle $\lambda(\E)$ on $S$ 
equipped with the Quillen metric $h^{\lambda(\E)}$. If the family 
$(\pi:\X\to S,\tilde\omega,\L)$ is effective, then by Fujiki and Schumacher
\cite{FS}, one can define the generalized Weil-Petersson metric 
$\widehat{h}_{WP}$ on $S$ and it follows that the first Chern form 
$c_1(\lambda(\E),h^{\lambda(\E)})$ of the determinant line bundle satisfies
$$c_1(\lambda(\E),h^{\lambda(\E)})=a_n\widehat\omega_{WP}$$
where $a_n=-2^{n+1}\kappa_2(n+1)!$. This result extends to the moduli space
${\mathfrak M}_{H,e}$ of extremal Hodge manifolds \cite{FS}.

\smallskip
\noindent
{\bf Acknowledgement.} The authors would like to thank Kai K\"ohler for
some very helpful comments and remarks.

\section[Preliminaries]{Preliminaries}
\setcounter{equation}{0}

Let $(X,g)$ be a compact K\"ahler manifold of complex dimension $n$. We denote
by
$\omega=\omega_g$  the K\"ahler form of $g$. Then the volume form of $g$
is given by   $dv_g=\omega^n_g/n!$.
We always choose the holomorphic connection on $X$. The
Riemann and Ricci curvature tensors will be denoted by $R$ and $r$,
respectively, and the scalar curvature by $s_g$. In terms of the Ricci tensor $r$, the Ricci form $\rho$
is defined by
$$\rho(\varphi,\psi)=r(J\varphi,\psi),$$
where $J$ denotes the complex structure of $X$. For a K\"ahler manifold,
$i\rho$ is the curvature of the canonical line bundle
$K=(\Lambda^{n}T^{1,0}X)^*$.
This has important implications for the scalar curvature $s_g=\Tr(r)$. It can
be calculated by the formula
\begin{equation}\label{1.1}
s\thinspace\omega^n=2n\;\rho\wedge\omega^{n-1}.
\end{equation}
Since $\rho/2\pi$ represents the first Chern class $c_1(X)=c_1(T^{1,0}X)
=-c_1(K)$, it follows that the average value 
\begin{equation}\label{1.2}
s_0=\frac{\int_X s_g dv_g}{\int_X dv_g}=4n\pi\frac{\int_X
c_1(\Omega)\wedge\omega^{n-1}}{\int_X\omega^n}
\end{equation}
of the scalar curvature is a topological invariant, depending only on the
K\"ahler class $[\omega]$.

\smallskip

\noindent
Let $\E\to X$ be a holomorphic vector bundle with hermitian metric
$h^\E$. Let 
\begin{equation*}
0\to \Lambda^{0,0}(X,\E)\stackrel{\overline\partial}{\longrightarrow}\Lambda^{0,1}(X,\E)\stackrel{\overline\partial}{\longrightarrow}
\cdots\stackrel{\overline\partial}{\longrightarrow}\Lambda^{0,n}(X,\E)\to
0
\end{equation*}
be the Dolbeault complex and let 
$\Delta^\E_{0,q}= \overline\partial\thinspace \overline\partial^* 
+\overline\partial^*\overline\partial$ be the corresponding Dolbeault-Laplace
operator acting in $\Lambda^{0,q}(X,\E)$. Let $P^\E_{0,q}$ be
the  projection on harmonic $(0,q)$-forms and for 
$\Re(s)>n/2$, let
\begin{equation*}
\zeta_q(s;\E)=\frac{1}{\Gamma(s)}\int^\infty_0t^{s-1}\mbox{Tr}
(e^{-t\Delta^\E_{0,q}}-P^\E_{0,q})dt
\end{equation*}
be the zeta function of $\Delta^\E_{0,q}$. It admits a meromorphic
continuation to $\C$, also denoted by $\zeta_q(s;\E)$,  which is regular at 
$s=0$. Then the zeta regularized determinant $\det\Delta^\E_{0,q}$ of
$\Delta^\E_{0,q}$ is defined
by
\begin{equation*}
\det\Delta^\E_{0,q}=\exp\left(-\frac{d}{ds}\zeta_q(s;\E)\big|_{s=0}\right).
\end{equation*}
\smallskip

\noindent
The Ray-Singer analytic torsion is defined to be 
\begin{equation*}
T_0(X,\E)=\left(\prod^n_{q=0} (\det\Delta^\E_{0,q})^{(-1)^{q+1} q}\right)^{1/2}
\end{equation*}
\cite{RS}. This number depends, of course, on the metrics $(g,h^\E)$
and if it is necessary to indicate this dependence we shall write
$T_0(X,\E,g,h^\E)$ for the analytic torsion. 
Let 
\begin{equation*}
\lambda(\E)=\bigotimes^n_{q=0} (\det H^q(X,\E))^{(-1)^{q+1}}
\end{equation*}
be the determinant line associated to the Dolbeault complex. Let
$\H^{0,q}(\E)$  be the space of $\E$-valued harmonic
$(0,q)$-forms.
For each $q$, there is a canonical isomorphism 
${\H}^{0,q}({\E})\stackrel{\cong}{\longrightarrow}H^q(X,\E)$. 
Thus, using the $L^2$-metric on 
${\H}^{0,q}(\E)$, we get a metric $\parallel\cdot\parallel_{L^{2}}$ on
$\lambda(\E)$. The Quillen metric on $\lambda(\E)$ is then defined by
\begin{equation}\label{1.3}
\parallel\cdot\parallel_Q=\parallel\cdot\parallel_{L^{2}} \cdot T_0(X,\E).
\end{equation}
We shall also consider virtual holomorphic bundles $\E=\sum_{k=1}^mn_k\E_k$,
where $n_k\in\Z$ and $\E_k\to X$ are holomorphic vector bundles. 
For such
$\E$ we set
\begin{equation*}
\lambda(\E)=\bigotimes_{k=1}^m\lambda(\E_k)^{n_k}.
\end{equation*}
Here for $n_k<0$, $\lambda(\E_k)^{n_k}:=(\lambda(\E_k)^*)^{-n_k}.$
If each $\E_k$ is equipped with a hermitian metric, we get a metric
$\parallel\cdot\parallel_Q$ on $\lambda(\E)$ which is the tensor product of
the induced Quillen metrics on $\lambda(\E_k)^{n_k}$.
\smallskip

\noindent
Let $v\in \lambda(\E)$, $v\neq 0$. Then we define
\begin{equation*}
{\mathcal Q}(X,\E,g,h):=\parallel v\parallel_{Q,g,h}.
\end{equation*}
If $X$ and $\E$ are fixed, we shall  write $\cQ(g,h)$ in place of
$\cQ(X,\E,g,h)$.

Let $\omega$ be a fixed K\"ahler form. Put
$$\K_{\omega}=\bigl\{g\mid g\text{ K\"ahler metric on } X\text{ with }[\omega_g]=[\omega]\bigr\}.
$$
Let 

\begin{equation}\label{1.4}
\H_{\omega}=\bigl\{\varphi\in C^\infty(X,\R)\mid \omega+i\partial\overline
\partial\varphi>0\bigr\}
\end{equation}
be the corresponding space of K\"ahler potentials. If $\varphi\in\H_{\omega}$,
then 
\begin{equation}\label{1.5}
\omega(\varphi):=\omega+i\partial\overline\partial\varphi
\end{equation}
is the 
K\"ahler form of a K\"ahler metric $g(\varphi)\in\K_\omega$ and it is well
known that the
natural map
\begin{equation}\label{1.6}
\varphi\in\H_\omega\to g(\varphi)\in\K_\omega
\end{equation}
is surjective. 
If $\varphi\in\H_\omega$
then for each $c\in\R$, one has $\varphi+c\in\H_\omega$ and
$\omega(\varphi)=\omega(\varphi+c)$. On the other hand, if
$\varphi,\psi\in\H_\omega$ are such that $\omega(\varphi)=\omega(\psi)$ then 
$\partial\overline\partial(\varphi-\psi)=0$ and therefore $\varphi-\psi$ is
constant. Let 
$$\H_{\omega}^0=\left\{\varphi\in\H_\omega\;\Big|\;
  \int_X\varphi\omega^n=0\right\}.$$
Then the restriction of the map (\ref{1.6}) to the subspace $\H^0_\omega$
  induces an isomorphism
\begin{equation}
\H_\omega^0\cong\K_\omega.
\end{equation}
Furthermore, for any $\varphi\in C^\infty(X)$ there exists $\epsilon>0$ such
that 
\begin{equation*}
\omega_u=\omega+iu\partial\overline\partial\varphi>0
\end{equation*}
for all $u\in\R,|u|<\epsilon$. The corresponding family $g_u,|u|<\epsilon$, of
K\"ahler metrics will be called the variation of $g$ in the 
$\varphi$-direction.

\section[The variational formala]{The variational 
formula}
\setcounter{equation}{0}
Given a smooth family $g_t$, $t\in \R$, of K\"ahler metrics, we will 
denote by
$*_t$ the Hodge star operator with respect to $g_t$ and we set
\begin{equation}\label{2.1}
U_t:= (g_t)^{-1} \frac{d}{dt}
(g_t),
\quad \alpha_t:=*^{-1}_t \frac{d}{dt}
(*_t).
\end{equation}

\noindent
We observe that if $g_t$ is the variation of $g$ in the $\varphi$-direction, 
then we have 
\begin{equation}\label{2.2}
\mbox{Tr}(U_t)= \sum_{\alpha,{\overline\beta}} g^{\alpha{\overline\beta}}_t
\frac{\partial^2\varphi}{\partial z_\alpha \partial{\overline z}_\beta}
=-\frac{1}{2}\Delta^{t}_{0,0}\varphi.
\end{equation}
The variation of ${\mathcal Q}(X,\E,g,h)$ has been computed by Bismut, Gillet
and Soul\' e
\cite[Theorem 1.22]{BGS3}. We recall their result. Let $g_t,t\in\R$,
be a smooth family of K\"ahler metrics on $X$ and $h_t,\; t\in\R$, a smooth
family of hermitian metrics on $\E$. Let $\Omega_t$ and $L^\E_t$ be the
curvature of the Hermitian holomorphic connection on $(T^{(0,1)} X,g_t)$ and
$(\E,h_t)$, respectively. Let Td be the Todd genus and $Td_j$ its $j$-th
component. The Todd genus is normalized as in \cite{BGS3}. Then
\begin{equation}\label{2.3}
\begin{split}
\frac{\partial}{\partial t} \log {\mathcal Q}(X,\E,g_t,h_t) &  =
\frac{1}{2}\left( \frac{1}{2\pi i}\right)^n\int_X\frac{\partial}{\partial b}
\bigg[ \mbox{Td}(-\Omega_t-b U_t) \\ &  \times
\mbox{Tr}
\left( \exp(-L^\E_t -b(h^\E_t)^{-1}
\frac{\partial h^\E_t}{\partial t})\right)\bigg]_{b=0}.
\end{split}
\end{equation}

In particular, if we assume that $\E$ is trivial, this formula is reduced to 

\begin{equation}\label{2.4}
\frac{\partial}{\partial t}\log {\mathcal Q}(X,g_t)=
\frac{1}{2}\left(\frac{1}{2\pi i}\right)^n \int_X
\frac{\partial}{\partial b}
\left[ \mbox{Td}_{n+1}(-\Omega_t-bU_t)\right]\Big|_{b=0}.
\end{equation}

Using the definition of ${\mathcal Q}(X,g)$ and (1.120) of \cite{BGS3},
it follows that
\begin{equation}\label{2.5}
\frac{\partial}{\partial t}\log{\mathcal Q}(X,g_t)=
\frac{\partial}{\partial t}\log T_0(X,g_t)-\frac{1}{2}
\sum^n_{q=0}(-1)^q\Tr(\alpha_t P^t_{0,q}).
\end{equation}

For the application that we have in mind we need a more explicit version of
the variational formula. Recall that the Todd genus can be expressed in
terms of the Chern classes. Therefore, we have to compute the corresponding
derivatives of the Chern classes which are also normalized as in \cite{BGS3}.

\begin{lem}\label{l2.1} 
For all $j\ge1$, we have
\begin{equation}\label{2.6}
\frac{\partial}{\partial b}\left( {c_j}(-\Omega-b
  U)\right)\Big|_{b=0}= \sum^{j-1}_{k=0}(-1)^{j+k}\Tr(\Omega^k U)
{c}_{j-k-1}(\Omega).
\end{equation}
\end{lem}

\begin{pf}
We proceed by induction on $j$. For $j=1$ we have
\begin{equation*}
\frac{\partial}{\partial b}\left( c_1(-\Omega-b
  U)\right)\Big|_{b=0} = -\mbox{Tr}(U)
\end{equation*}
proving (\ref{2.6}) in this case. Now suppose that (\ref{2.6}) holds for
$j$. We express the Chern classes by the power sums $s_k$. If we 
formally write
\begin{equation*}
\sum^n_{j=0} c_j(\Omega)x^j=\prod^n_{j=1}(1+\gamma_j x),
\end{equation*}
then $s_k$ is the $k$-th elementary symmetric function in
$\gamma_1,\ldots,\gamma_n$. Using Newton's formula, we obtain
\begin{equation*}
 c_{j+1}(-\Omega-b U)=\frac{1}{j+1}\sum^{j+1}_{k=1} 
(-1)^{k+1} (c_{j+1-k} s_k)(-\Omega-b U).
\end{equation*}
Since $s_k(R)=\mbox{Tr}(R^k)$, we get
\begin{equation*}
\frac{\partial}{\partial b}(s_k(-\Omega -b U))\mid_{b=0} 
=k(-1)^k\mbox{Tr}(\Omega^{k-1}U).
\end{equation*}
Using the induction hypothesis, we deduce

\begin{equation*}
\begin{split}
 \frac{\partial}{\partial b}\bigl( &{c}_{j+1}(-\Omega-b U)\bigr)
\Big|_{b=0}\\
&=\frac{1}{j+1}\sum^{j+1}_{k=1}
\left\{ \sum^{j-k}_{l=0}(-1)^{j+l}
\mbox{Tr}(\Omega^lU) c_{j-k-l}(\Omega) s_k(-\Omega)\right.\\
&\hskip130pt -k{ c}_{j+1-k}(-\Omega)\mbox{Tr}(\Omega^{k-1} U)\Biggr\}\\
& =\frac{1}{j+1} \sum^j_{l=0} \mbox{Tr}(\Omega^l U)
  \left\{ \sum^{j-l}_{k=1} (-1)^{l+j+k} {
    c}_{j-k-l}(\Omega)s_k(\Omega)\right.\\
&\hskip143pt -(-1)^{j+l}(l+1){ c}_{j-l}(\Omega)\Biggr\}\\
&=   \sum^j_{l=0} (-1)^{j+1+l}\mbox{Tr}(\Omega^l U)
{c}_{j-l}(\Omega).
\end{split}
\end{equation*} 
\end{pf}

Expressing the Todd genus in terms of the
Chern classes and using Lemma \ref{l2.1}, it follows from (\ref{2.4}) that
there exist $a_{j_{0}\cdots j_{n+1}}\in\C$, $j_0,\ldots,j_{n+1}\in\N$, 
such that 
for all $\varphi\in C^\infty(X):$
\begin{equation}\label{2.7}
\begin{split}
\frac{\partial}{\partial t}\log{\mathcal Q}(X,g_t)
&=\sum_{j_0+\cdots+j_{n+1}=n}
a_{j_{0}\cdots j_{n+1}}\\
&\quad\times\int_X \Tr(\Omega^{j_{n+1}}_t U_t)
{c}_{j_{0}}(\Omega_t)\cdots {c}_{j_{n}}(\Omega_t).
\end{split}
\end{equation}

The coefficients $a_{ j_{0}\cdots j_{n+1}}\in\C$, $j_0,\ldots,j_{n+1}\in N$, 
can
 be computed explicitly. For example, if $n=2$, then formula
(\ref{2.7}) takes the form
\begin{equation}\label{2.8}
\begin{split}
\frac{\partial}{\partial t}\log{\mathcal Q}(X,g_t)&
=\frac{1}{48}\int_X\Tr(\Omega_tU_t)c_1(\Omega_t)\\
&-\frac{1}{48}\int_X\Tr(U_t)(c_2(\Omega_t)+c_1(\Omega_t)^2).
\end{split}
\end{equation}

Let $\varphi\in C^\infty(X)$ and let $g_u$, $|u|<\epsilon$, 
 be the variation of $g$ in
the direction of $\varphi$. Put
\begin{equation}\label{2.9}
\frac{\delta}{\delta\varphi}\log{\mathcal Q}(X,g)=\frac{\partial}{\partial u}
\log {\mathcal Q}(X,g_u)\Big|_{u=0}.
\end{equation}
If $n=2$, it folows from (\ref{2.2}) and (\ref{2.8}) that
\begin{equation}\label{2.10}
\begin{split}
\frac{\delta}{\delta\varphi}\log{\mathcal Q}(X,g)&
=\frac{1}{48}\int_X\Tr(\Omega U_0)c_1(\Omega)\\
&+\frac{1}{96}\int_X\Delta\varphi(c_2(\Omega)+c_1(\Omega)^2).
\end{split}
\end{equation}

\section[Riemann surfaces]{Riemann surfaces}
\setcounter{equation}{0}

In this section we consider the case $n=1$. Then $X$ is a compact oriented
surface without boundary. First observe that

\begin{equation*}
\alpha(1)=-\frac{\Delta\varphi}{2}\; \mbox{and}\; \alpha({d\overline z})=0.
\end{equation*}
Hence, it follows from (\ref{2.5}) that

\begin{equation*}
\frac{\delta}{\delta\varphi}\log{\mathcal
  Q}(X,g)=\frac{\delta}{\delta\varphi}\log T_0(X,g).
\end{equation*}
Furthermore, $\overline\partial: \wedge^{0,0}(X)\to\wedge^{0,1}(X)$ is an
isomorphism on nonzero eigenspaces. This implies that 
$\det \Delta_{0,1}=\det\Delta_{0,0}$. Let $\Delta_g=d^* d$ denote the
Laplacian on functions. Then $\Delta_g=2\Delta_{0,0}$ and therefore we get
\begin{equation}\label{3.1}
2\log T_0(X,g)=\log \det\Delta_g+(\log 2)\left( 1-\frac{\chi(X)}{6}\right).
\end{equation}
Thus there exists $c>0$ such that
\begin{equation*}
T_0(X,g)=c(\det\Delta_g)^{1/2}.
\end{equation*}
Next we describe the space ${\mathcal K}_\omega$ for a given K\"ahler class
$[\omega ]\in H^{1,1}(X)$. Since $H^{1,1}(X)\cong\R,\; [\omega]$ is unique up to
multiplication by $\R^+$. Let $g$ be the metric with K\"ahler form $\omega$
and let $\Delta=\Delta_g$. Let
$\psi\in C^\infty(X)$. Then $\omega+i\partial\overline\partial \psi>0$ holds
if and only if $1-\Delta\psi >0$. Thus the space  $\H_\omega$, which is 
 defined by (\ref{1.4}), is given by
\begin{equation}\label{3.2}
{\mathcal H}_\omega \cong \{ \psi\in C^\infty(X)\mid 1-\Delta\psi>0\}.
\end{equation}
 Then by using $g(\psi)=(1-\Delta\psi)g$,  we may regard $T_0(X,g)$ as a
functional on ${\mathcal H}_\omega$ and by (\ref{3.1}) we have
\begin{equation*}
 2 \log T_0(X,\psi)=\log \det\Delta_{(1-\Delta\psi)g}+C,\quad
  \psi\in{\mathcal H}_\omega.
\end{equation*}
Let
\begin{equation*}
\mbox{Conf}(g)=\{e^{2\varphi}g\mid\varphi\in C^\infty(X)\}
\end{equation*}
be the conformal equivalence class of metrics on $X$ which contains $g$. Let
\begin{equation*}
\mbox{Conf}_0(g)=\{{\widetilde g}\in\mbox{Conf}(g)\mid \mbox{Area}(X,{\widetilde
  g})=\mbox{Area}(X,g)\}.
\end{equation*}

\begin{lem}\label{l3.1}
Let $g$ be a K\"ahler metric on $X$ with K\"ahler form $\omega$. Then 
\begin{equation*}
{\mathcal K}_\omega=\Cf_0(g).
\end{equation*}
\end{lem}

\begin{pf} Let ${\widetilde g}\in{\mathcal K}_\omega$. By (\ref{3.2}) there
  exists a unique $\psi\in C^\infty(X)$ with $\int\psi\;dv_g=0$ such that
  ${\widetilde g}=(1-\Delta\psi)g$. 
Since $1-\Delta\psi>0$, there exists $\varphi\in
  C^\infty(X)$ such that $e^{2\varphi}=1-\Delta\psi$. Moreover
\begin{equation*}
\int_X e^{2\varphi}dv_g=\int_Xdv_g-\int_X \Delta\psi dv_g=\int_Xdv_g.
\end{equation*}
Hence there exists a unique $\varphi\in C^\infty(X)$ such that
${\widetilde g}=e^{2\varphi}g\in{\mbox{Conf}}_0(g)$. On the other hand,
  let ${\widetilde g}\in\mbox{Conf}_0(g)$. Then ${\widetilde g}=e^{2\varphi}g$
  and $\int_X(e^{2\varphi}-1)dv_g=0$.
Therefore, there exists a unique $\psi\in C^\infty(X)$ with $\int\psi\;dv_g=0$
such that
  $e^{2\varphi}=1-\Delta\psi>0$. 
Hence
\begin{equation*}
{\widetilde g}=(1-\Delta\psi)g\in{\mathcal K}_\omega.
\end{equation*}
\end{pf}

\noindent
The functional 
\begin{equation}\label{3.3} 
{\widetilde g}\in\mbox{Conf}_0(g)\longmapsto\log\det\Delta_{ {\widetilde g}}\in\R
\end{equation}
has been studied by Osgood, Philipps and Sarnak \cite{OPS1}. By (\ref{3.1})
and Lemma 3.1, we see that, up to a constant, (\ref{3.3}) coincides
with ${\widetilde g}\in{\mathcal K}_\omega\longmapsto \log{\mathcal
  Q}(X,{\widetilde g})$. We shall now derive the main  result of \cite{OPS1} 
by our approach.
\smallskip

\noindent
Let $g_t\in\K_\omega$, $t\in\R$, be a smooth family of metrics. Then there
exists a smooth family $\psi_t\in C^\infty(X)$ such that $g_t=g(\psi_t)$.
Using Lemma 2.1, the variational formula (\ref{2.4}) gives
\begin{equation}\label{3.4}
\frac{\partial}{\partial t}\log{\mathcal Q}(X,g_t)=-\frac{1}{24}\int_X
c_1(\Omega_t)\Delta_t\dot \psi_t=\frac{1}{12\pi}\int_X K_t\Delta_t\dot \psi_t
\;dv_t, 
\end{equation}
where $K_t$ is the Gauss curvature of $g_t$. We note that the last equality
follows from 
$$c_1(\Omega_t)=-\frac{2K_t}{\pi}dv_t.$$ In particular, we get
\begin{equation}\label{3.4a}
\frac{\delta}{\delta\psi}\log{\mathcal Q}(X,g)=-\frac{1}{24}\int_X
c_1(\Omega)\Delta \psi=\frac{1}{12\pi}\int_X K\Delta \psi\;
dv. 
\end{equation}
\smallskip

\begin{remark}
This formula should be compared with the corresponding
variational 
formula  (1.12) in
\cite{OPS1}. Given $\psi\in C^\infty(X)$, there exists $\epsilon>0$ such that
$1-u\Delta \psi>0$ for $|u|< \epsilon$. Then by Lemma \ref{l3.1}, there exists
$\phi_u\in 
C^\infty(X)$ such that 
\begin{equation}\label{3.5}
1-u\Delta\psi=e^{2\phi_u},\quad |u|<\epsilon.
\end{equation}
The family $g_u=e^{2\phi_u}g,\quad |u|<\epsilon$, of metrics is a
deformation of $g$ in $\mbox{Conf}_0(g)$. By (\ref{3.5}), we have
\begin{equation*}
\phi_0=0,\quad {\dot\phi}_0=-\frac{\Delta\psi}{2}.
\end{equation*}
Then with respect to the variation defined by $g_u$, in formula (1.12)
of \cite{OPS1} we have $\delta\varphi={\dot\phi}_0$, $\varphi=0$ and
$\delta \log A=0$.
Hence, up to a constant, (\ref{3.4a}) corresponds to (1.12) of \cite{OPS1}. 
\end{remark}

By (\ref{3.4a}), a metric ${\widetilde g}$ is critical for 
$\log{\mathcal Q}(X,g)$
if and only if its Gauss curvature $K$ is constant. To investigate 
 the critical point, we use the evolution of a Riemannian metric $g$ on $X$
by Hamilton's  Ricci flow equation with $s_0$ as in (\ref{1.2})
\begin{equation}\label{3.6}
\frac{\partial g_{ij}}{\partial t}=(s_0-s_g)g_{ij}.
\end{equation}
Note that by the Gauss-Bonnet formula
$$s_0=\frac{4\pi\chi(X)}{A},$$
where $A=\text{Area}(X,g)$.
In \cite{Ha}, Hamilton proved that if $s_0\le0$ or if $s_g>0$, then for any
initial data, 
(\ref{3.6}) has a unique solution $g(t)$ which exists for all time and
as $t\to\infty$,  converges to a metric of constant curvature.
The case $s_0>0$ was completed by Chow \cite{Ch} who proved that for $s_0>0$,
the scalar curvature 
$s(t)$ of the solution $g(t)$ of (\ref{3.6}) becomes  positive in finite
time. 
 Now if $g(0)$ is hermitian
then $g(t)$ is hermitian for all time. Furthermore, the Ricci flow is area
preserving \cite[p.238]{Ha}. Therefore, the Ricci flow preserves the space 
$\K_\omega$. Hence, there exists $u(t)\in C^\infty(X)$ which is
a smooth function of $t\in\R^+$ such that, with respect to a local holomorphic
parameter $z$,
$$g(t)=g(0)+\frac{\partial^2u(t)}{\partial z\partial\overline{z}}dz\otimes
 d\overline{z}.$$
Let $g(t)=h(t)dz\otimes d\overline{z}$. 
Then (\ref{3.6}) implies
\begin{equation}\label{3.7}
\Delta_t \dot u=-h^{-1}\frac{\partial^2\dot u}{\partial z\partial\overline{z}}
=(s_g-s_0).
\end{equation}
Combined with (\ref{3.4}) we obtain
\begin{equation}\label{3.8}
\begin{split}
\frac{\partial}{\partial t}\log\cQ(g(t))&=\frac{1}{24\pi}\int_Xs(t)\Delta_t\dot u(t)\;dv_t\\
&=\frac{1}{24\pi}\int_X(s(t)-s_0)\Delta_t\dot u(t)\;dv_t\\
&=\frac{1}{24\pi}\int_X(s(t)-s_0)^2\;dv_t\ge0,
\end{split}
\end{equation}
where $s(t)$ denotes the scalar curvature of the metric $g(t)$. If $s_0\le0$,
Hamilton \cite{Ha} proved that there exist $\varepsilon>0$ and $C>0$ such that
\begin{equation}\label{3.10}
|s(t)-s_0|\le Ce^{-\varepsilon t},\quad t\in\R^+.
\end{equation}
If $s_0>0$, he  established the same estimate under the additional assumption
that  $s_g>0$. By Chow's results \cite{Ch} it follows that (\ref{3.10}) holds
in the case $s_0>0$ too. Hence we can integrate (\ref{3.8}) and we get
$$\log\cQ(g(\infty))-\log\cQ(g(0))=\frac{1}{24\pi}
\int_0^\infty\int_X(s(t)-s_0)^2\;dv_t\;dt\ge0.$$

Furthermore, if $s(0)$ is not constant then $(s(t)-s_0)^2>0$ for $t$ in a
neighborhood of $0$ and therefore,
the above inequality is strict. Thus we proved the following theorem.
\begin{thm}
$\cQ$ has a unique maximum on $\K_\omega$ which is attained at the metric
$g^*\in\K_\omega$ of constant curvature.
\end{thm} 
This result was first proved by Osgood, Phillips and Sarnak in \cite{OPS1}
by a different method.

\section[Constant holomorphic curvature]{Manifolds with constant 
holomorphic curvature}
\setcounter{equation}{0}
 
Let $(X,g)$ be a K\"ahler manifold and let $J$ denote the complex structure of
 $X$. Recall that $(X,g)$ is said to have
constant holomorphic curvature $H\in\R$, if for any $z\in X$ and any $\xi\in
T_z X$ with $\langle \xi,\xi\rangle=1$, the holomorphic curvature
\begin{equation*}
H_z(\xi)=\langle R_z(\xi,J\xi)\xi,J\xi\rangle
\end{equation*}
is equal to $H$. We note that any metric of constant holomorphic curvature is
K\"ahler-Einstein \cite[(6.12)]{Go}. 
Furthermore, compact manifolds with constant holomorphic curvature $H$ can be
classified as follows \cite{Be}: If $H>0$, then $X$ is isomorphic to
$\C P^n$ equipped with the rescaled Fubini-Study metric; if $H=0$, then $X$ is
isomorphic to a complex torus $\Gamma\backslash \C^n$ with the flat metric; if
$H<0$ then $X$ is isomorphic to a compact quotient $\Gamma\backslash D^n$ of
the complex unit ball $D^n\subset \C^n$, equipped with a Bergmann metric, by a
discrete group $\Gamma$ of isometries. Thus manifolds with constant
holomorphic curvature may be regarded as higher-dimensional analogues of
Riemann surfaces. So, if we think of generalizing the results of Osgood,
Philipps and Sarnak \cite{OPS1} to higher dimensions, these are the most
obvious candidates for 
being critical points of the functional $\log{\mathcal Q}(X,g)$. 
\smallskip

\noindent
We wish to apply the variational formula (\ref{2.4}). For this purpose we need
several lemmas.

\begin{lem}\label{l4.1}
Suppose that $g$ is a K\"ahler metric on $X$ with constant holomorphic
curvature $H$, and let $\Omega$ denote the curvature form of $g$. Then
\begin{enumerate}
\item With respect to local holomorphic coordinates $(z_1,\ldots,z_n)$,
  $\Omega$ is given by
\begin{equation*}
\Omega^j_k=-iH\delta_{jk}\omega+\frac{H}{2}\sum^n_{m=1}
g_{k{\overline m}} dz_j\wedge d\overline z_m.
\end{equation*}
\item The Chern classes of $(X,g)$ are represented by
\begin{equation*}
c_j(\Omega)= \binom{n+1}{j}
\left( -\frac{H}{2\pi}\omega\right)^j.
\end{equation*}
\end{enumerate}
\end{lem}

\begin{pf}
1.  The claimed formula for $\Omega^j_k$ is an immediate consequence of the
  following formula for the Riemann curvature tensor $R$ on a manifold with
  constant holomorphic curvature $H$ \cite[(6.1.1)]{Go}:
\begin{equation*}
R_{j{\overline k}l{\overline m}}=\frac{H}{2}\left(
g_{j{\overline k}} g_{l{\overline m}}+g_{j{\overline m}}g_{l{\overline
    k}}\right) .
\end{equation*}

2.  We recall that a representative of the $j$th Chern class is given in
  terms of the curvature form $\Omega$ by means of the formula
\begin{equation*}
c_j(\Omega)=\left(\frac{1}{2\pi i}\right)^j\frac{1}{j!}
\sum_{\begin{Sb}
k_1,\ldots,k_j,l_1,\ldots,l_j\\ \in\{1,\ldots,n\}
\end{Sb}}
\delta^{l_1\cdots l_j}_{k_1\cdots
  k_j}\Omega^{k_1}_{l_1}\wedge\cdots\wedge\Omega^{k_j}_{l_j} 
\end{equation*}
\cite[(6.12.1)]{Go}. Inserting the expression for $\Omega^j_k$, given in 1.,
a straightforward computation leads to 
\begin{equation*}
c_j(\Omega)= \left(\frac{H}{2\pi i}\right)^j\sum^j_{p=0} \binom{n-j+p}{p}
\left( -i\omega\right)^j= \binom{n+1}{j}\left( -\frac{H}{2\pi}\omega\right)^j.
\end{equation*}

\end{pf}

 Let $U=U_0$, where $U_t$ is defined  by (\ref{2.1}).

\begin{lem}\label{l4.2}
Let $(X,g)$ be the K\"ahler manifold with constant holomorphic curvature
$H$. Let $\varphi \in C^\infty(X)$ and let 
$\omega_u=\omega+iu\partial{\overline\partial}\varphi,\quad
|u|<\epsilon$. Then for all $j\ge1$,
\begin{equation*}
\mbox{Tr} (\Omega^j U)=-\frac{\Delta \varphi}{2}(-iH\omega)^j+
\partial{\overline\partial}\left( \frac{H}{2}\varphi(-iH\omega)^{j-1}\right).
\end{equation*}
\end{lem}

\begin{pf}
Choose local holomorphic coordinates $(z_1,\ldots,z_n)$. Then we claim that for
all $j\ge1$:
\begin{equation}\label{4.2}
\begin{split}
(\Omega^j U)^p_q= &  \sum^n_{k=1}\left\{ (-i H\omega)^jg^{p{\overline k}}
\varphi_{q{\overline k}}\right.\\
& + (-iH\omega)^{j-1}\frac{H}{2}\left. \varphi_{q{\overline k}}dz_p\wedge
d{\overline z}_k\right\},
\end{split}
\end{equation}
where 
$\varphi_{q{\overline k}}=\partial^2/\partial z_q \partial
{\overline z}_k\varphi$.
To prove (\ref{4.2}), we proceed by induction on $j$. 
If $j=1$, then using Lemma \ref{l4.1}, we get
\begin{equation*}
\begin{split}
(\Omega U)^p_q   &   =\sum^n_{l=1}\Omega^p_l U^l_q=
 \sum_{l,k}\left( -iH\delta_{pl}\omega+\frac{H}{2}\sum^n_{m=1}
g_{l{\overline m}}dz_p\wedge dz_{\overline m}\right)
g^{l{\overline k}}\varphi_{q{\overline k}}\\
& =\sum^n_{k=1}\left( -iH\omega g^{p{\overline k}}\varphi_{q{\overline k}}+
\frac{H}{2} \varphi_{q{\overline k}}dz_p\wedge{\overline z}_k\right).
\end{split}
\end{equation*}
Suppose that (\ref{4.2}) holds for $j$. Then using Lemma \ref{l4.1} combined
with the induction hypothesis, we obtain
\begin{equation*}
\begin{split}
(\Omega^{j+1}U)^p_q&  =\sum^n_{l=1}\Omega^p_l(\Omega^jU)^l_q
= \sum_{kl}\left( -iH\delta_{pl}\omega+\frac{H}{2}\sum^n_{m=1}
g_{l{\overline m}}dz_p \wedge d{\overline z}_m\right)\\
&\wedge\left( (-iH\omega)^jg^{l{\overline k}}\varphi_{q{\overline
      k}}+(-iH\omega)^{j-1}\frac{H}{2} \varphi_{q{\overline k}}
dz_l\wedge d{\overline z}_k\right)\\
&=\sum^n_{k=1}\left\{ (-iH\omega)^{j+1}g^{p{\overline k}} \varphi_{q{\overline
      k}}+(-iH\omega)^j\frac{H}{2}\varphi_{q{\overline k}}dz_p\wedge
  d{\overline z}_k\right\}.
\end{split}
\end{equation*}
This proves (\ref{4.2}). Thus in local holomorphic coordinates, we get
\begin{equation*}
\begin{split}
\mbox{Tr}(\Omega^jU)&  =(-iH\omega)^j\sum^n_{p,k=1}g^{p{\overline k}}
\frac{\partial^2\varphi}{\partial z_p\partial{\overline z}_k}+
\frac{H}{2}(-iH\omega)^{j-1} \sum^n_{p,k=1}
\frac{\partial^2\varphi}{\partial z_p\partial{\overline z}_k} 
dz_p\wedge d{\overline z}_k\\
& = - \frac{\Delta\varphi}{2}(-i H\omega)^j+\frac{H}{2}(-iH\omega)^{j-1}
\partial{\overline\partial}\varphi.
\end{split}
\end{equation*}
Since $\omega$ is closed, this implies the lemma.

\end{pf} 

\begin{thm}\label{th4.3}
Let $g_{KE}\in\K_\omega$ be a metric of constant holomorphic curvature. 
Then $g_{KE}$ is
a critical point of $\cQ\colon\K_\omega\to\R$.
\end{thm}

\begin{pf}  
By (\ref{1.6}) it suffices to prove that 
\begin{equation*}
\frac{\delta}{\delta \varphi }\log{\mathcal Q}(X,g_{KE})=0
\end{equation*}
for all $\varphi\in C^\infty(X)$. Let $\varphi\in C^\infty(X)$ and let $g_u$,
$|u|<\varepsilon$, be the corresponding variation of $g_{KE}$ in the
$\varphi$-direction. We use the variational formula (\ref{2.7}) and apply 
 Lemma \ref{l4.1} and Lemma \ref{l4.2} to the integrant. Then it
follows that there exist $A,B\in \C$ such that
\begin{equation*}
\begin{split}
\frac{\delta}{\delta\varphi}  \log{\mathcal Q}(X,g_{KE})&  =
\int_X\left\{
  A\Delta\varphi\omega^n+B\partial{\overline\partial}(\varphi\omega^{n-1})
\right\}\\
& =
A\int_X\Delta\varphi\omega^n+B\int_X\partial
{\overline\partial}(\varphi\omega^{n-1}).
\end{split}
\end{equation*}
Hence, we get
$$\frac{\delta}{\delta\varphi}  \log{\mathcal Q}(X,g_{KE})=0.$$
This concludes the proof.
\end{pf}

\noindent
In the case of a complex torus, Theorem 0.1 can be restated in terms of the
analytic torsion. Namely, we have

\begin{cor}\label{c4.4} Let $X=\Gamma\backslash \C^n$ be a complex torus and
  let $g_{KE}$ be the flat metric on $X$. Then 
\begin{equation*}
\frac{\delta}{\delta\varphi} T_0(X,g_{KE})=0
\end{equation*}
for all $\varphi\in C^\infty(X)$.
\end{cor}

\begin{pf}
Let $\varphi\in C^\infty(X)$ and let $g_u,\; |u|<\epsilon$, be the variation
of $g_{KE}$ in the $\varphi$-direction. Let $\alpha=\alpha_0$ denote the operator
(\ref{2.1}) at $u=0$ where $*_u$ is the Hodge star operator with respect to
$g_u$. Let $P_{0,q}, \; 0\le q\le n$,  denote the harmonic projections with
respect to $g_{KE}$. By Theorem \ref{th4.3}
 and (\ref{2.5}) it is sufficient to prove that
\begin{equation*}
\sum^n_{q=0 } (-1)^q Tr(\alpha P_{0,q})=0
\end{equation*}
for all $\varphi\in C^\infty(X)$. Let $(z_1,\ldots,z_n)$ denote the standard
coordinates on $X$ induced from $\C^n$. Recall that an orthonormal basis for
the space $\H^{0,q}(X)$ of harmonic $(0,q)$-forms is given by
\begin{equation*}
\left\{ (\mbox{vol}(X))^{-1/2} d{\overline z}^B \mid B\subset \{
  1,\ldots,n\},\  |B|=q\right\}.
\end{equation*}
Using the definition of the Hodge star operator and the assumption that
$(g_{KE})_{i{\overline j}}=\delta_{ij},$ it follows that

\begin{equation*}
\begin{split}
\overline{ Tr(\alpha P_{0,q})} & = c\sum_{|B|=q}\int 
\frac{\partial}{\partial u}
\left\{ \det \left(( g^{a{\overline b}}_u)_{a,b\in
    B}\right)\det\left((g_u)_{i{\overline j}}\right)\right\}\Big|_{u=0} \;dv_0\\
& = c \sum_{|B|=q} \int_X\left\{ \sum_{b\in B} {\dot g}^{b{\overline b}}-
    \frac{\Delta \varphi}{2}\right\}\;dv_0  \\
&= c_1\int_X\Delta\varphi\;dv_0=0.
\end{split}
\end{equation*}
\end{pf}

\begin{remark}
If $\K_\omega$ contains a metric $g_{KE}$ of constant holomorphic curvature,
then we have seen that $g_{KE}$ is a critical point of  
$\log{\mathcal Q}(X,\tilde g)$ for $\tilde g\in{\mathcal K}_{\omega}$.
We do not know, however, wether or not $g_{KE}$ is the unique critical
point of 
$\log{\mathcal Q}$ on ${\mathcal K}_{\omega}$. 
Also, we do not know anything about the nature of the critical point. It would 
be interesting to see wether or not
$g_{KE}$ is  a maximum or a minimum of $\log{\mathcal Q}(X,\tilde g)$.
\end{remark}

\section[K3 surfaces]{K3 surfaces}
\setcounter{equation}{0}
Another class of examples are $K3$ surfaces. Recall that a  $K3$ surface  is a
compact, connected complex analytic surface $X$ that is regular, meaning
$h^1(X,{\mathcal O})=0$, and its canonical bundle $K=\Lambda^{2,0}(T^*X)$ is
trivial. By Siu \cite{Si}, every $K3$ surface admits a K\"ahler
metric. Furthermore, by Yau \cite{Y}, every K\"ahler class $[\omega]$ of $X$
contains a unique K\"ahler-Einstein metric $g_{KE}$. This metric should be 
a natural candidate for a critical point of $\cQ$ and we shall now investigate
under what conditions this is the case.
\begin{prop}\label{prop5.1}
Let $X$ be a $K3$ surface and let $[\omega]$ be a K\"ahler class of $X$. Then
the unique K\"ahler-Einstein metric $g_{KE}\in\K_\omega$ is a critical point
of $\cQ:\K_\omega\to\R$ if and only if $\Delta c_2(\Omega_{KE})=0$, where 
$\Omega_{KE}$ is the curvature form of $g_{KE}$. 
\end{prop}
\begin{pf}
Since $c_1(X)=0$, the K\"ahler-Einstein metric $g_{KE}$ is Ricci
flat. Therefore we have $c_1(\Omega_{KE})=0$. Then it follows from the
variational formula (\ref{2.10}) that
$$\frac{\delta}{\delta\varphi}\log\cQ(g_{KE})=\frac{1}{96}\int_X\Delta\varphi
c_2(\Omega_{KE})=\frac{1}{96}\int_X\varphi\Delta c_2(\Omega_{KE})$$
for all $\varphi\in C^\infty(X)$. Hence $g_{KE}$ is a critical point of $\cQ$
if and only if $\Delta c_2(\Omega_{KE})=0$.
\end{pf}

\noindent
The condition $\Delta c_2(\Omega_{KE})=0$ can be reformulated in terms of the
curvature tensor.

\begin{lem}\label{lem5.2}
Let $g$ be a Ricci flat metric on $X$. Let $R$ be the curvature tensor of
$g$. Then we have
$$c_2(\Omega)=\frac{1}{4\pi^2} |R|^2\omega^2,$$
where $| R|$ is the pointwise norm of $R$.
\end{lem}
\begin{pf}
 Let $W$ denote
the Weyl curvature tensor of $g$. Since $g$ is Ricci flat, it follows from 
\cite[(1.116)]{Be} that $R=W$.
 Furthermore, let $L:\Lambda^{p,q}(T^*X)\to\Lambda^{p+1,q+1}(T^*X)$ be the
 operator 
$L(\eta)=\omega\wedge \eta$ and let $\Lambda$ be its adjoint. Then it follows
from equation (2.80) in \cite{Be} that
$$\Lambda^2(c_2(\Omega))=\frac{1}{4\pi^2}|W|^2.$$
 This implies the lemma.
\end{pf}
Combining Proposition \ref{prop5.1} and Lemma \ref{lem5.2}, we obtain
\begin{cor}\label{cor5.3}
The K\"ahler-Einstein metric $g_{KE}\in\K_\omega$ is a critical point of $\cQ$
if and only if $|R(x)|$ is constant.
\end{cor}

\smallskip
\noindent
The condition that $|R(x)|$ is constant is very stringent and we do not know if
there exists any $K3$ surface which admits a K\"ahler-Einstein metric
satisfying this condition. However, by a result of S. Kobayashi \cite{Ko}, one
knows that on
certain $K3$ surfaces there exist K\"ahler-Einstein metrics with curvature
concentrated near some divisor. These are examples of  K\"ahler-Einstein 
metrics on K3 surfaces such that  $|R|$ is not constant.


\section[Algebraic manifolds]{Algebraic manifolds}
\setcounter{equation}{0}
In this section we introduce a twisted version of our functional  $\cQ$ which
has a simpler variational formula. 
For this purpose  we now assume that $X$ is a projective algebraic
manifold. Note 
that by the Kodaira embedding theorem \cite{GH} any compact complex manifold
with definite Chern class is a projective algebraic manifold. 
\smallskip

\noindent
A polarization of $X$ is the choice of an ample line bundle $L\to X$ up to
isomorphism. The pair $(X,L)$ is called a polarized algebraic manifold. Since
$L$ is ample, a certain power $L^{\otimes k}$, $k\in\N$, defines an embedding
$i_L: X\hookrightarrow \C P^n$ \cite{GH}. Let $g$ be the pullback of the
Fubini-Study metric on $\C P^n$ and let $[\omega]$ be the K\"ahler class
determined by $g$. Then $[\omega]$ is a rational multiple of $c_1(L)$. So we
can normalize $g$ such that $[\omega]=c_1(L)$. This implies that 
$[\omega]\in H^2(X,\Z)\cap H^{1,1}(X,\R)$.  
We fix a hermitian metric $h$ on $L$ such that the curvature $\Theta_h$
satisfies $\Theta_h=-2\pi i\omega$. Such metrics $h$ exist \cite[pp.163,191]{GH}
and, up to multiplication by a positive real number, $h$ is uniquely
determined by $\omega$. Given $\varphi\in C^\infty(X,\R)$, put 
\begin{equation*}
h(\varphi)=e^{-2\pi\varphi}h.
\end{equation*}
Then $h(\varphi)$ is a hermitian metric on $L$ whose curvature is given by
\begin{equation*}
\Theta_{h(\varphi)}=\Theta_h+2\pi\partial{\overline\partial}\varphi=-2\pi
i(\omega +i\partial{\overline{\partial}}\varphi).
\end{equation*}
Thus, if $\varphi\in{\mathcal H}_\omega$, we have
\begin{equation}\label{6.1}
\Theta_{h(\varphi)}=-2\pi i\omega(\varphi).
\end{equation}
In this section we consider the analytic torsion with coefficients in a
certain virtual vector bundle $\E$ associated to $L$. For its definition we
introduce the following numbers
\begin{equation*}
\kappa_1=\int_Xc_1(\Omega)\wedge \omega^{n-1} \;\;\mbox{and}\;\;
\kappa_2=\int_X\omega^n.
\end{equation*}
Since $[\omega]\in H^2(X,\Z)$, it follows that $\kappa_1$ and $\kappa_2$ are
integers. Then by (\ref{1.2}), we have
$$s_0=4n\pi\frac{\kappa_1}{\kappa_2}.$$
Now we define the virtual vector bundle $\E$ over $X$ to be
\begin{equation*}
\E=\bigoplus^{4\kappa_2(n+1)} (L-L^{-1})^{\otimes n}\oplus
\bigoplus^{-\kappa_1n}(L-L^{-1})^{\otimes(n+1)}.
\end{equation*}
If $\widetilde h$ is any hermitian metric on $L$, we denote by 
${\widetilde  h}^{\E}$ the induced hermitian metric on $\E$. 
Using this notation we
  introduce the following functional on $\H_\omega$:
\begin{equation}\label{6.3}
\cQ_L(X,\varphi)={\mathcal Q}(X,\E,g(\varphi),h(\varphi))^{(-1)^n},\quad 
 \varphi\in\H_\omega.
\end{equation}
Since $X$ is fixed, we set $\cQ_L(\varphi)=\cQ_L(X,\varphi)$.
To begin with we compute the variation of $\cQ_L(\varphi)$.

\begin{thm}\label{th6.1} 
Let $\varphi_u,\;u\in[a,b]$, be a smooth family of functions in 
${\mathcal H}_\omega$. Let
$\omega_u=\omega+i\partial\overline\partial\varphi_u$ and let $s(u)$ be the
scalar curvature of the metric $g_u=g(\varphi_u)$. Then
\begin{equation*}
\frac{\partial}{\partial u}\log \cQ_L(\varphi_u)=c_n\int_X{\dot\varphi}_u
(s(u)-s_0)\omega^n_u,
\end{equation*}
where 
$$c_n=\kappa_2(n+1)2^{n-1}.$$
\end{thm}

\begin{pf} To compute the variation, we apply formula (\ref{2.3}). Let 
 $F^{\E}_u$ denote the curvature form of the metric $h^{\E}_u$, induced
on $\E$ by $h(\varphi_u)$, and let  $V_u=(h^\E_u)^{-1} \frac{d}{du}(h^{\E}_u)$.
Then by (\ref{2.3}) we have

\begin{equation}\label{6.4}
\begin{split}
\frac{\partial}{\partial u}\log \cQ_L(\varphi_u)&  =\frac{(-1)^n}{2}
\left(\frac{1}{2\pi i}\right)^n\int_X\left[ \frac{\partial}{\partial
  b}\Td(-\Omega_u-bU_u)\Big|_{b=0}\ch(\E,h^\E_u)\right.\\
&   +\left.\Td(-\Omega_u)\frac{\partial}{\partial b}
\ch (-F^\E_u-b V_u)\Big|_{b=0}\right].
\end{split}
\end{equation}

\noindent
We start with the computation of the first term. Using that 
$\Theta_u=\Theta_{h_{u}}=-2\pi i \omega_u$, it follows that

\begin{equation*}
\begin{split}
  \ch(\E,h^\E_u)&=4\kappa_2(n+1)(e^{2\pi i\omega_u}-e^{-2\pi i\omega_u})^n-
\kappa_1 n(e^{2\pi i\omega_u} -e^{-2\pi i\omega_u})^{n+1}\\
 &   =4\kappa_2(n+1)(4\pi i\omega_u)^n.
\end{split}
\end{equation*}
Hence, by Lemma \ref{l2.1} we obtain
\begin{equation}\label{6.5}
\begin{split}
  \int_X\biggl[ \frac{\partial}{\partial b}\Td(&-\Omega_u-bU_u)\Big|_{b=0}
\ch(\E,h^\E_u)\biggr]^{(n)}\\
&=2\kappa_2(n+1)(4\pi i)^n\int_X\frac{\partial}{\partial
  b}c_1(-\Omega_u-bU_u)\Big|_{b=0}\omega_u^n\\
& =\kappa_2(n+1)(4\pi i)^n\int_X(\Delta_u\varphi_u)\omega^n_u=0.
\end{split}
\end{equation}
As for the second term, we have 
\begin{equation*}
\begin{split}
\frac{\partial}{\partial b}\ch \bigl( &-F^\E_u-bV_u\bigr) \Big|_{b=0}\\
&=(n+1)\left\{ 4n\kappa_2 \ch (L_u^{-1}-L_u)^{n-1}-n\kappa_1 \ch(L_u^{-1}-L_u)^n\right\}\\
& \hskip20pt\cdot \frac{\partial}{\partial b}\left\{ \ch(2\pi i\omega_u+b 2\pi
  {\dot\varphi}_u)-\ch(-2\pi i\omega_u-b2\pi\dot\varphi_u)\right\}\Big|_{b=0}\\
& =2\pi (n+1)\left\{ 4n\kappa_2(e^{-2\pi i\omega_u}
-e^{2\pi i\omega_u})^{n-1}\right.\\
&\hskip43pt\left.-n\kappa_1(e^{-2\pi i\omega_u}-e^{2\pi i\omega_u})^n\right\}\dot\varphi_u
 (e^{2\pi i\omega_u}+ e^{-2\pi i\omega_u})\\
& =4\pi (n+1)\kappa_2(-4\pi i)^n\left\{-\frac{n}{\pi i}
\omega^{n-1}_u-n\frac{\kappa_1}{\kappa_2}\omega^n_u\right\} \dot\varphi_u.\\
\end{split}
\end{equation*}

\noindent
By (\ref{1.1}), we have $c_1(\Omega_u)=\rho/2\pi$. Hence, using (\ref{1.2}),
we get  
\begin{equation*}
\begin{split}
\Biggl[ \Td(-&\Omega_u)\frac{\partial}{\partial
    b}\ch(-F^\E_u-bV_u)\Big|_{b=0}\Biggr]^{(n)}\\
& = 4\pi (n+1)\kappa_2(-4\pi i)^n\left(nc_1(\Omega_u)\wedge
    \omega^{n-1}_u-n\frac{\kappa_1}{\kappa_2}\omega^n_u\right)\dot\varphi_u\\
&=2^{n}(n+1)\kappa_2(2\pi i)^n(-1)^n\dot\varphi_u(s(u)-s_0)\omega_u^n.
\end{split}
\end{equation*}

\noindent
Together with (\ref{6.4}) and (\ref{6.5}) this leads to
\begin{equation*}
\frac{\partial}{\partial u}\log \cQ_L(\varphi_u)=c_n\int_X \dot\varphi_u (s(u)-s_0)
\omega_u^n.
\end{equation*}

\end{pf}

\begin{remark} 
There are other possible choices for $\E$ which give essentially
the same variational formula. Tian \cite{Ti1} used in a different context a
virtual bundle of the form
$$\bigoplus^a\left((K_X^{-1}-K_X)\otimes(L-L^{-1})^{n}\right)\oplus
\bigoplus^b(L-L^{-1})^{n+1}.$$
Up to a constant, it gives the same variational formula.
\end{remark}

\noindent
Let $\varphi_1,\varphi_2\in{\mathcal H}_\omega$. Let 
$\{\varphi_u\mid a\le u\le b\} $ be a piecewise smooth path in ${\mathcal
  H}_\omega$ such that $\varphi_a=\varphi_1$ and $\varphi_b=\varphi_2$. Then
it follows from Theorem \ref{th6.1} that 
\begin{equation}\label{6.6}
\log \left( \frac{\cQ_L(\varphi_1)}{\cQ_L(\varphi_2)}\right)=-c_n\int^b_a
\left\{ \int_X \dot\varphi_u (s(u)-s_0)\omega^n_u\right\} \; dt.
\end{equation}

\noindent
The right hand side of (\ref{6.6})  is equal to $c_n$ times the functional
$M(\varphi_1,\varphi_2)$ introduced by Mabuchi \cite[(2.2.2)]{Ma1}. It is
defined for any compact K\"ahler manifold $X$. 
For a projective algebraic manifold $X$, (\ref{6.6}) implies that
$M(\varphi_1,\varphi_2)$ is independent of the path 
$\{ \varphi_u\mid a\le u\le b\}$ in ${\mathcal H}_\omega$ 
connecting $\varphi_1$ and $\varphi_2$. This was proved by Mabuchi
\cite[Theorem 2.4]{Ma1} in general, using different methods.
 
\smallskip
\noindent
Let $\varphi\in{\mathcal H}_\omega$ and
$c\in\R$. Set $\varphi_u=\varphi+uc$, $u\in[0,1]$. Then by (\ref{6.6}) and 
(\ref{1.2}) we get 
\begin{equation*}
\begin{split}
  \log\left( \frac{\cQ_L(\varphi+c)}{\cQ_L(\varphi)}\right)&=c_nc\int^1_0\left\{
    \int_X(s(u)-s_0)\omega^n_u\right\}\;du\\
& = \frac{c_nc}{n!} \int^1_0 \left\{ \int_Xs(u)dv_g-s_0\mbox{Vol}(X)\right\}
    \;du=0.
\end{split}
\end{equation*}

\noindent
Hence for all $\varphi\in{\mathcal H}_\omega$ and $c\in\R$ we have
\begin{equation*}
\cQ_L(\varphi+c)=\cQ_L(\varphi).
\end{equation*}
Thus by (\ref{1.6}), $\cQ_L:{\mathcal H}_\omega\to \R$ factors through ${\mathcal
  K}_\omega$. The induced functional on ${\mathcal K}_\omega$ will be denoted
  by 
  the same letter:
\begin{equation*}
\tilde g\in{\mathcal K}_{\omega}\longmapsto \cQ_L(\tilde g)\in\R.
\end{equation*}
By (\ref{6.6}), the map $\mu\colon{\mathcal K}_{\omega_{0}}\to\R$ defined by 
\begin{equation}\label{6.7}
\mu(\tilde g)= c_n^{-1}\log \frac{\cQ_L(g)}{\cQ_L(\tilde g)}
\end{equation}
coincides with the $K$-energy map defined by 
Mabuchi \cite[Section3]{Ma1}. 

\noindent
Let $g_u\in\K_\omega$,
$|u|< \epsilon$, be a smooth path. Then by (\ref{1.6}) there exists a smooth
path $\varphi_u\in\H_\omega$ such that $g_u=g(\varphi_u)$. It follows 
from Theorem \ref{6.1} that
\begin{equation*}
\frac{\partial}{\partial u}\log\cQ_L(g_u)\Big|_{u=0}=c_n
\int_X\dot\varphi_0(s_{g_0}-s_0)\omega_0^n.
\end{equation*}
In particular, if $g_u$, $|u|<\epsilon$, is the variation of $g$ in the
direction of some $\varphi\in C^\infty(X)$, then this formula implies 
\begin{equation*}
 \frac{\delta}{\delta\varphi}\log\cQ_L(g)=c_n \int_X\varphi(s_g-s_0)\omega^n.
\end{equation*}
This proves the following theorem.

\begin{thm}\label{th6.2}
A K\"ahler metric $g_0\in{\mathcal K}_{\omega}$ is a critical point of
${\mathcal Q}_L\colon\K_{\omega}\to\R$ if and only if the scalar curvature
$s_{g_0}$ 
of $g_0$ is constant. In this case $s_{g_0}=s_0$.
\end{thm}

\smallskip
\noindent
This holds for the K-energy in general \cite{Ma1}.

\smallskip
\noindent
In  \cite{Ca}, Calabi introduced the notion of extremal K\"ahler metrics where
a 
K\"ahler metric $\tilde g\in{\mathcal K}_{\omega}$ is called extremal if 
$\tilde g$ is
a critical point of the functional
\begin{equation*}
S(\widetilde g)=\int_X s^2_{\widetilde g}\; dv_{\widetilde g},\quad
{\widetilde g}\in {\mathcal K}_{\omega}.
\end{equation*}
Any metric $\tilde g$ of constant scalar curvature is extremal, but as shown by
Calabi \cite{Ca}, the converse is not true, i.e., there are extremal K\"ahler
metrics with nonconstant scalar curvature. However, if $X$ is non-uniruled,
then the extremal K\"ahler  metrics are precisely the metrics with constant
scalar 
curvature. 
\smallskip

\noindent
Now we collect some further information about the functional $\cQ_L$.
By Theorem \ref{th6.2}, the critical points of $\cQ_L$ are the  metrics of
constant scalar curvature. If $\K_\omega$ contains a K\"ahler-Einstein metric,
then 
the subset $\K_E\subset \K_\omega$ of all K\"ahler-Einstein metrics coincides
with 
the set of critical points of $\cQ_L$. In particular, if $c_1(X)\le0$, then by
\cite{A}, \cite{Y}, $\cQ_L$ has critical points. On the other hand, the Futaki
invariant \cite{Fu}, \cite{Ca2} obstructs the existence of constant scalar
curvature metrics in certain K\"ahler classes. Therefore, critical points do
not always exist.

The K-energy map $\mu\colon\K_\omega\to\R$ has been studied by Mabuchi and
Bando 
\cite{Ma1}, \cite{Ma2}, \cite{B}, \cite{BM}, mainly in connection with the
Futaki obstruction to the existence of K\"ahler-Einstein metrics on compact
complex manifolds with $c_1(X)>0$. Mabuchi \cite{Ma2} defined  a natural
Riemannian structure on $\K_\omega$, i.e., $\K_\omega $ can be equipped
canonically with  
the structure of an
infinite-dimensional Riemannian manifold. The tangent space 
 $T_{\tilde\omega}\K_\omega$ of 
$\K_\omega$ at $\tilde\omega\in\K_\omega$ can be identified with
$$\left\{\eta\in C^\infty(X,\R)\mid \int_X\eta\tilde\omega^n=0\right\}$$
and the Riemannian structure is given by
$$\langle\eta_1,\eta_2\rangle=\frac{1}{\text{Vol(X)}}
\int_X\eta_1\eta_2\frac{\tilde\omega^n}{n!},\quad\eta_1,\eta_2\in 
T_{\tilde\omega}\K_\omega.$$
One of the main results of  \cite[Theorem 5.3]{Ma2} states 
 that $\mu\colon\K_\omega\to\R$ is a
convex function, meaning that $\text{Hess}\:\mu$ is positive semidefinite.
Therefore, by (\ref{6.7}) the same holds for $\log\cQ_L$. This resembles the
situation in the Riemann surface case.

If $\K_\omega$ contains a K\"ahler-Einstein metric, then by Theorem \ref{0.1},
$\cQ_L$ is bounded from above. In general, we do not know anything about
boundedness of $\cQ_L$.

Finally, we note that by Theorem \ref{6.1} the gradient flow of $\log\cQ_L$ is
given by
$$\frac{\partial\varphi }{\partial t}(t)=s(t)-s_0,$$
where $s(t)$ is the scalar curvature of the metric $g(t)=g(\varphi(t))$.
Using that
$$s=2\sum_{kl}g^{k\overline l}r_{k\overline l}\quad\text{and}\quad
r_{k\overline l}=-\frac{\partial^2}{\partial z_k\partial\overline z_l}
\log\det(g_{\alpha\overline \beta}),$$
we get the following fourth order nonlinear parabolic equation
\begin{equation}\label{6.8}
\frac{\partial\varphi}{\partial t}=-2\sum_{k,l}g(t)^{k\overline l}\frac{\partial^2}{\partial
  z_k\partial\overline{z_l}}\left(\log\det\left(g_{\alpha\overline \beta}
+\frac{\partial^2\varphi}{\partial
  z_\alpha\partial\overline{z_\beta}}\right)\right)-s_0.
\end{equation} 
Differentiating this equation with respect to $t$, we get
$$\frac{\partial}{\partial t}\left(\frac{\partial\varphi}{\partial
  t}\right)=-2\Delta_t^2\left(\frac{\partial\varphi}{\partial
  t}\right)+\sum_{k,l}r(t)^{k\overline l}\frac{\partial^2}{\partial
  z_k\partial\overline{z}_l}\left(\frac{\partial \varphi}{\partial
  t}\right),$$
where $r(t)$ is the Ricci tensor of the metric $g(t)$. So, from standard
  theory we know that the solution of the initial value problem for
  (\ref{6.8})
exists for short time. The crucial question is to prove existence for all
time. If $\varphi(t)$ is any solution of (\ref{6.8}), then by Theorem
\ref{6.1} we get
$$\frac{\partial}{\partial t}\log\cQ_L(g(\varphi(t))=
c_n\int_X(s(t)-s_0)^2\omega_t^n\ge0.$$
Thus  $\log\cQ_L$ is increasing along the gradient flow.

\section[K\"ahler-Einstein metrics]{K\"ahler-Einstein metrics}
\setcounter{equation}{0}

\noindent
By Theorem \ref{th6.2}, the critical points of $\cQ_L$ are exactly the metrics
of 
constant scalar curvature. In general, it seems to be difficult to determine
the 
nature of the critical points. Much more can be said if $[\omega]$ contains
a K\"ahler-Einstein metric $\omega_{KE}$. Recall that $\omega$ is said to be
K\"ahler-Einstein if the Ricci form $\rho_\omega$ is proportional to the
K\"ahler form
\begin{equation}\label{7.1}
\rho_\omega=\lambda \omega
\end{equation}
for some $\lambda\in\R$. Since $2\pi c_1(X)=[\rho_\omega]$, the first Chern
class has to satisfy either one of the following conditions
\begin{equation*}
c_1(X)=0\quad \mbox{or}\quad c_1(X)>0\quad\mbox{or}\quad c_1(X)<0.
\end{equation*}
If $c_1(X)=0$, it follows from Yau \cite{Y} that each K\"ahler class
$[\omega]$ contains a unique K\"ahler-Einstein metric $\omega_{KE}$. If
$c_1(X)<0$, i.e., if the canonical line bundle $K_X$ is ample, then it was
shown by T. Aubin \cite{A} and S.-T. Yau \cite{Y}, that up to multiplication by
a positive scalar, there exists a unique K\"ahler-Einstein metric $g_{KE}$ on
$X$ which can be normalized such that
$2\pi c_1(X)=-[\omega_{KE}]$. In the case $c_1(X)>0$ there are obstructions
for the 
existence of K\"ahler-Einstein metrics \cite{Fu}, \cite{Ti2} and the existence
problem has not been  completely settled yet. See \cite{Bo} for a review of
the results.
First we make the following elementary observation.

\begin{lem}\label{l7.1} Let $(X,g)$ be a compact K\"ahler manifold with 
 definite
  or vanishing first Chern class. Suppose that there is a K\"ahler-Einstein
  metric $g_{KE}$ on $X$ with $[\omega_{KE}]=[\omega_g]$. Then $g$ has
  constant scalar
  curvature iff $g$ is K\"ahler-Einstein.
\end{lem}
\begin{pf} It follows from (\ref{1.1}) that a K\"ahler-Einstein metric has
  constant scalar curvature.
Now assume that $g$ has constant scalar curvature $s_g$. Then $s_g=s_0$. Let 
$g_{KE}$ be a K\"ahler-Einstein metric with $[\omega_{KE}]=[\omega]$. By
(\ref{1.1}), we have 
\begin{equation*}
\rho_{KE}=\frac{s_0}{2n}\omega_{KE}.
\end{equation*}
Since $[\rho_{KE}]=[\rho]$, we get $\frac{s_0}{2n}[\omega]=[\rho]$. Hence there
exists 
$\varphi \in C^\infty(X,\R)$ such that
\begin{equation*}
\rho=\frac{s_0}{2n}\omega+i\partial{\overline\partial}\varphi.
\end{equation*}
Taking the trace of both sides of this equation yields
\begin{equation*}
\frac{s_0}{2}-\frac{1}{2}\Delta\varphi=\text{tr}(r)=\frac{s_g}{2}
=\frac{s_0}{2}.
\end{equation*}
This implies that $\varphi$ is constant. Therefore, we get 
$\rho=\frac{s_0}{2n}\omega$.
\end{pf}

\begin{cor}\label{cor7.2}
Let $\K_{KE}\subset\K_\omega$ denote the set of all K\"ahler-Einstein metrics 
contained in $\K_\omega$. Suppose that $\K_{KE}\not=\emptyset$. Then 
$\K_{KE}$ is precisely the set of critical points of $\cQ_L$.
\end{cor}

\smallskip
\noindent
{\bf Proof of Theorem 0.1.}
First assume that $c_1(X)\le 0$. Since by our assumption, the set of 
K\"ahler-Einstein metrics which are contained in $\K_\omega$ is nonempty, it 
follows from \cite{A} and \cite{Y} that $\K_\omega$ contains a unique 
K\"ahler-Einstein metric $g_{KE}$, and if $c_1(X)<0$, then
$\frac{1}{2\pi}[\omega_{KE}]$ represents $-c_1(X)$. Now it follows from 
 Corollary \ref{cor7.2} that $g_{KE}$  is the unique critical point of
$\cQ_L:\K_{\omega}\to\R$. To prove that $g_{KE}$ is the absolut maximum
of $\cQ_L$, we employ the complex
analogue of the Ricci flow.  
We normalize $g_{KE}$ such that
\begin{equation*}
\rho_{KE}=\frac{s_0}{2n}\omega_{KE}.
\end{equation*}
 Let $g\in\K_{\omega}$. By H.-D. Cao \cite{C}, the evolution equation
\begin{equation}\label{7.2}
\frac{\partial {\widetilde g}_{i{\overline j}}}{\partial t}
=-{\widetilde r}_{i{\overline j}}+\frac{s_0}{2n}{\widetilde g}_{i{\overline
    j}},\quad
{\widetilde g}_{i{\overline j}}(0)=g_{i{\overline j}},\quad i,j=1,\ldots,n,
\end{equation}
has a unique solution ${\widetilde g}(t)$ which exists for all $t\ge 0$ and
satisfies
\begin{equation}\label{7.3}
\lim_{t\to \infty} {\widetilde g}(t)=g_{KE}
\end{equation}
in the $C^\infty$-topology. Moreover, there exists $u\in C^\infty(\R^+\times
X,\R)$ with
\begin{equation}\label{7.4}
{\widetilde g}_{i{\overline j}}(t,z)=g_{i{\overline j}}(z)+
\frac{\partial^2u(t,z)}{\partial z_i\partial{\overline z}_j}.
\end{equation}
Let $\Delta_t$ be the Laplace operator with respect to ${\widetilde g}(t)$. By
(\ref{7.2}) and (\ref{7.4}) we obtain
\begin{equation*}
\begin{split}
\frac{1}{2}\Delta_t({\dot u}(t))& =-\sum_{i,j}{\widetilde g}^{i{\overline
    j}}(t) \frac{\partial^2{\dot u}(t)}{\partial z_j\partial {\overline
    z}_j}=-\sum_{i,j}{\widetilde g}^{i{\overline j}}(t)\frac{\partial
    {\widetilde g }_{i{\overline j}}}{\partial t}(t)\\
&= \sum_{i,j}{\widetilde g}^{i{\overline j}}(t) \left( \widetilde 
r_{i{\overline j}}(t)
-\frac{s_0}{2n}
    {\widetilde g}_{i{\overline j}}(t)\right)=\frac{{\widetilde s}(t)}{2}-
\frac{s_0}{2}.
\end{split}
\end{equation*}
Thus $u(t)$ satisfies
\begin{equation}\label{7.5}
\Delta_t(\dot u(t))=\widetilde s(t)-s_0.
\end{equation}
Together with Theorem \ref{th6.1} we get

\begin{equation}\label{7.6}
\begin{split}
\frac{\partial}{\partial t}\log \cQ({\widetilde g}_t) & = c_n\int_X{\dot
  u}(t)({\widetilde s}(t)-s_0)\omega^n_t\\
& = c_n\int_X{\dot u}(t)\Delta_t({\dot u}(t))\omega^n_t\\
&= \frac{c_n}{n!}\int_X|\nabla \dot u(t)|_t^2\;dv_t.
\end{split}
\end{equation}
Hence, by the definition of $c_n$, we obtain
\begin{equation}\label{7.7}
\frac{\partial}{\partial t}\log\cQ(\widetilde g_t)\ge0.
\end{equation}
Suppose that $\widetilde s(0)$ is not constant. Then by (\ref{7.5}) there
exists 
$\varepsilon>0$ such that for all $t\le\varepsilon$, $\dot u(t,z)$ is not
constant  as a function of $z$. By (\ref{7.6}) it follows that
\begin{equation}\label{7.8}
\frac{\partial}{\partial t}\log\cQ(\widetilde g_t)>0\quad\text{for all}\quad
t\le\varepsilon.
\end{equation}

By \cite[Proposition 2.2]{C}, ${\dot u}(t)$ converges to a constant in the
$C^\infty$-topology as $t\to\infty$. Also ${\widetilde g}_{i{\overline j}}(t)$
  converges to $(g_{KE})_{i{\overline j}}$ in the $C^\infty$-topology as
  $t\to\infty$ \cite[Main Theorem]{C}. This implies that there
  exists $C_1>0$ such that
\begin{equation}\label{7.9}
\sup_{z\in X}| \Delta_t({\dot u}(t,z))|\le C
\end{equation}
for all $t\ge 0$.
Put
\begin{equation*}
a(t)=\frac{1}{\mbox{Vol}(X)} \int_X{\dot u}(t) dv_t, \; t\ge 0.
\end{equation*}
By (2.11) and (2.18) of \cite{C} there exist $C_2,C_3>0$ such
that 
\begin{equation}\label{7.10}
\int_X|{\dot u}_t-a(t)|dv_t\le C_2 e^{-C_3t}.
\end{equation}
Using (\ref{7.9}) and (\ref{7.10}), we obtain
\begin{equation}\label{7.11}
\left| \int_X{\dot u}_t\Delta({\dot u}_t)dv_t\right|=
\left| \int_X({\dot u}_t-a(t))\Delta_t({\dot u}_t)dv_t\right|
\le C_4 e^{-C_3t}.
\end{equation}
Since for $t\to\infty$, $\widetilde g_{i\overline j}(t)$ converges to
$(g_{KE})_{i\overline j}$ in the $C^\infty$-topology, it follows that
\begin{equation}\label{7.13}
\log\cQ_L(g_{KE})=\lim_{t\to\infty} \log \cQ_L({\widetilde g}_t).
\end{equation}
Combing (\ref{7.6}), (\ref{7.11}) and (\ref{7.13}), it follows that for 
all $t>0$, we have
\begin{equation*}
\log \cQ(g_{KE})-\log \cQ({\widetilde g}_t)=c_n\int^\infty_t
  \frac{\partial}{\partial 
  v}\log \cQ({\widetilde g}_v)dv.
\end{equation*}
Hence, if $g\not=g_{KE}$, (\ref{7.7}) together with (\ref{7.8}) imply that
\begin{equation*}
\log \cQ_L(g_{KE})> \log \cQ_L(g).
\end{equation*}
This completes the proof of the first part of the theorem.

\smallskip
\noindent
Now assume that $c_1(X)>0$. Let $\K_{E}\subset\K_0$ be the subset of
K\"ahler-Einstein metrics. By our assumption $\K_{E}\not=\emptyset$. Then it
follows from \cite[Theorem 1]{B} that the K-energy map
$\mu\colon\K_\omega\to\R$ 
takes its absolute minimum on $\K_{E}$. Hence, by (\ref{6.7}), $\cQ_L$ 
attains its maximum on $\K_{E}$.

\section[Moduli spaces]{Moduli spaces and Quillen metric}
\setcounter{equation}{0}
So far, we considered the Quillen norm $\cQ_L$ of a fixed vector
$v\in\lambda(\E)$ as a function on $K_\omega$. In this section we investigate
the 
behaviour of the Quillen metric with respect to variations of the complex
structure. In \cite{FS}, Fujiki and Schumacher defined the moduli space of
extremal K\"ahler metrics on a fixed $C^\infty$ manifold. 

First we consider the local problem. Let $(\pi:\X\to S,\widetilde\omega)$ be a
metrically polarized family of compact K\"ahler manifolds \cite[Definition
3.2]{FS} over a reduced complex space $S$ and let $p\colon F\to \X$ be a
holomorphic hermitian  vector bundle with metric $h^F$. Let $\lambda(F)$ be the
holomorphic line bundle on $S$ associated to $(\det R\pi_*\F)^{-1}$, where
$\F$ is the locally free sheaf corresponding to $F$. Given $s\in S$, let
$X_s=\pi^{-1}(s)$ and $F_s=F| X_s$. We recall that there exists
a canonical isomorphism
\begin{equation}\label{8.1}
\lambda(F)_s\cong\lambda(F_s)=\bigotimes_{q\ge0}\left(\det
  H^q(X_s,F_s)\right)^{(-1)^{q+1}}
\end{equation}
\cite{BGS3}.
For each $s\in S$, $\tilde \omega$ defines a K\"ahler metric $g_s$ on $X_s$
and we get a family $\tilde g=\{g_s\}_{s\in S}$ of K\"ahler metrics. Let 
$h_s^F$ be the hermitian metric on $F_s$ induced by $h^F$. Let
$h^{\lambda(F_s)}$ be the Quillen metric on $\lambda(F_s)$ associated to
$(g_s,h^F_s)$. Using the isomorphism (\ref{8.1}), 
$\{h^{\lambda(F_s)}\}_{s\in S}$  defines a metric 
$h^{\lambda(F)}$ on $\lambda(F)$. It is proved in \cite{BGS3} that the metric
$h^{\lambda(F)}$ is   smooth. This metric is called  the {\it Quillen metric}
 on
$\lambda(F)$ associated to $(\tilde g,h^F)$.

If $S$ is a complex manifold, then the curvature of the
holomorphic connection on $\lambda(F)$ was computed by \cite[Theorem
1.27]{BGS3}. Namely the first Chern form is given by
\begin{equation}\label{8.2}
c_1(\lambda(F),h^{\lambda(F)})=-\left[\int_{\X/S}\Td(\X/S,\widetilde g)
\ch(F,h^F)\right]^{(2)}
\end{equation}
where $\Td(\X/S,\widetilde g)$ is the Todd class of the relative tangent
bundle. 
This formula was extended by Fujiki and Schumacher \cite[Theorem 10.1]{FS} to
the  case where $S$ is a reduced complex space.

Now let $(\pi:\X\to S,\widetilde\omega,\L)$ be a metrically polarized family
of compact Hodge manifolds \cite[Definition 3.8]{FS}. Thus $(\pi:\X\to
S,\L)$ is a polarized family of Hodge manifolds and $(\pi:\X\to
S,\widetilde\omega)$ is a metrically polarized family of K\"ahler manifolds
such that for each $s\in S$, $\omega_s$ represents $c_1(L_s)$ on $X_s$.
Furthermore, there exists a hermitian metric  $h$ on $\L$ such that the
restriction of $c_1(\L,h)$ to each fibre $X_s$ represents $\omega_s$
\cite[Proposition 3.10]{FS}.  Any such
metric is called admissible. 

Let
$$\kappa_1=\int_{X_s}c_1(X_s)\wedge\omega_s^{n-1}\quad\text{and}\quad\kappa_2=\int_{X_s}\omega_s^n.$$
These are integers and therefore, they  are independent of $s\in S$. Using
 $\L$, we introduce
the family version of the virtual holomorphic bundle used in Section 6. 
Let $n$ be the fibre dimension and set
$$\E=\bigoplus^{4(n+1)\kappa_2}(\L-\L^{-1})^n\oplus\bigoplus^{-n\kappa_1}(\L-\L^{-1})^{n+1}.$$
Let $h$ be an admissible hermitian metric on $\L$ for 
$(\pi,\widetilde\omega)$. Let $h^\E$ be the induced hermitian metric on
$\E$. Then the 
Chern character of $(\E,h^\E)$ is given by
\begin{equation*}
\begin{split}
\text{ch}(\E,h^\E)=4(n+1)\kappa_2\Biggl\{2^n&c_1(\L,h)^n
+\frac{n2^{n-1}}{3}c_1(\L,h)^{n+2}+\cdots\Biggr\}\\
&-n\kappa_1\left\{2^{n+1}c_1(\L,h)^{n+1}+\cdots\right\}.
\end{split}
\end{equation*}
Hence by (\ref{8.2}) we get the following expression for the first Chern form
of
$(\lambda(\E),h^{\lambda(\E)})$:
\begin{equation}\label{8.3}
\begin{split}
c_1(\lambda(\E),h^{\lambda(\E)})=-\kappa_2(n+1)&2^{n+1}\Biggl(\int_{\X/S}
c_1(\L,h)^nc_1(\X/S,\widetilde g)\\
&-\frac{n}{n+1}\frac{\kappa_1}{\kappa_2}\int_{\X/S}c_1(\L,h)^{n+1}\Biggr).
\end{split}
\end{equation}
Assume that the family $(\pi:\X\to S,\widetilde\omega,\L)$ is effective which
means that the Kodaira-Spencer map associated to $(\pi,\L)$ is injective 
\cite[Definition 4.1]{FS}. Then
by \cite[Definition 7.2]{FS}, there exists a generalized Weil-Petersson metric 
$\widehat h_{WP}=\{\widehat h_s\}_{s\in S}$ on $S$. Let $\widehat\omega_{WP}$
denote the corresponding Weil-Petersson form. Using Theorem 7.8 of \cite{FS}
together with (\ref{8.3}),
one gets the following theorem which is analogous to Theorem 10.3 of \cite{FS}.

\begin{thm}\label{th8.1} Let $(\pi:\X\to S,\widetilde\omega,\L)$ be an effective 
 family of extremal
  compact Hodge manifolds of constant scalar curvature with $S$ being
  connected, and let $h$ be an 
  admissible hermitian metric on $\L$ for $(\pi,\widetilde\omega)$. Then the
  first Chern form of the determinant line bundle
  $(\lambda(\E),h^{\lambda(\E)} )$ is given by
$$c_1(\lambda(\E),h^{\lambda(\E)} )=a_n\widehat\omega_{WP}$$
where $a_n=-2^{n+1}\kappa_2(n+1)!$.
\end{thm}

If we recall the isomorphism (\ref{8.1}) and the construction of the Quillen
metric on $\lambda(\E)$, we get an expression of the curvature in terms of the
analytic torsion. For each $s\in S$, the fibre $X_s$ is equipped with an
extremal Hodge metric $g_s$ and $\L_s$ is equipped with a hermitian metric
$h_s$ such that $\Theta_{h_s}=-{2\pi i}\omega_s$. Let $h^{\E_s}$ be the
associated hermitian metric in the virtual bundle 
$$\E_s=\bigoplus^{4(n+1)\kappa_2}(\L_s-\L_s^{-1})^n\oplus\bigoplus^{-n\kappa_1}(\L_s-\L_s^{-1})^{n+1}$$
and put $T_0(X_s,\E_s)=T_0(X_s,\E_s,g_s,h^{\E_s})$. Furthermore, let
$\parallel\cdot\parallel_{\lambda(\E_s)}$ denote the metric induced by
$(g_s,h^{\E_s})$ on the determinant line
$$\lambda(\E_s)=\bigotimes_{q\ge0}\left(\det
  H^q(X_s,\E_s)\right)^{(-1)^{q+1}}.$$ 
Now let 
 $\phi:S^\prime\to\lambda(\E)$ be a local holomorphic section of $\lambda(\E)$
  over some
open subset $S^\prime\subset S$. Using the definition of the Quillen metric
  in $\lambda(\E)$, it follows from Theorem \ref{th8.1} that
\begin{equation}\label{8.4}
\partial_s\overline\partial_s\log\left(\parallel\phi(s)\parallel^2_{\lambda(\E_s)}T_0(X_s,\E_s)\right)=a_n\widehat{\omega}_{WP}.
\end{equation}
In other words,
$\log\left(\parallel\phi(s)\parallel^2_{\lambda(\E_s)}T_0(X_s,\E_s)\right)$ is
a potential for the generalized Weil-Petersson metric on $S$.

\smallskip
\noindent
This construction can be globalized  \cite[\S 11]{FS}. Fix a compact
connected $C^\infty$ manifold $X$ and an integral class $\alpha\in H^2(X,\R)$.
Let ${\mathfrak M}_{H,e}$ be the moduli space of extremal Hodge manifolds with
underlying $C^\infty$ structure $(X,\alpha)$. Then ${\mathfrak M}_{H,e}$ is 
a complex $V$-manifold and the generalized Weil-Petersson metric is a $V$-form
$\omega_{WP}$ on ${\mathfrak M}_{H,e}$. As shown in \cite{FS}, the local
determinant line bundles can be patched together to a holomorphic line bundle
$F$ on ${\mathfrak M}_{H,e}$ and the Quillen metric on the local bundles
patches together to a $V$-metric $h^F$ on $F$. The global version of Theorem
\ref{th8.1} implies that  the Chern $V$-form  of $(F,h^F)$ satisfies
\begin{equation}\label{8.5}
c_1(F,h^F)=a\omega_{WP}
\end{equation}
for some integer $a\not=0$. In particular, this implies that the cohomology 
class
$[\omega_{WP}]$ is an integral cohomology class. Hence, any compact analytic
subspace of ${\mathfrak M}_{H,e}$ is a projective algebraic manifold.

We note that this is  well known for Riemann surfaces. Let $X$ be a compact
Riemann surface of 
genus $g\ge2$. Then there exists a discrete cocompact subgroup $\Gamma\subset
\text{SL}(2,\R)$ such that $X=\Gamma\backslash H$, where $H$ is the upper
half-plane. Let $g$ be the Riemannian metric on $X$ which is induced by the
Poincar\'e metric on $H$ and let
$\Delta=\overline\partial^*\overline\partial$ be the corresponding Laplace
operator. Fix a marking of $X$ and let $\tau$ be the corresponding period
matrix of $X$. We regard both $\Delta$ and $\tau$ locally as functions on the
moduli space ${\mathfrak M}_g$ of compact Riemann surfaces of genus $g$. Then
one has 
$$\partial_z\overline\partial_z\log\left(\frac{\det\Delta_z}
{\det\tau_z}\right)=\frac{i}{6\pi}\omega_{WP},$$
where $\omega_{WP}$ is the usual Weil-Petersson metric on ${\mathfrak
  M}_g$. This is  
essentially (\ref{8.5}). Furthermore, Wolpert \cite{Wo} has shown that
$\omega_{WP}$ extends to a closed form with singularities on the
campactification $\overline{\mathfrak M}_g$ of the moduli space and in the sense of
currents,  
$\frac{1}{\pi^2}\omega_{WP}$ is the Chern form of a continuous metric $h_{WP}$
on a certain line bundle $\lambda_{WP}$. This gives a projective embedding of
$\overline{\mathfrak M}_g$.

\end{document}